\documentclass{amsart}  

\usepackage{amsmath, latexsym, amssymb,stmaryrd}
\usepackage{color}
\usepackage{amsmath, amsfonts, amssymb, hyperref}
\usepackage{verbatim}
\usepackage[all]{xy}

\definecolor{gr75}{gray}{0.75}
\newcommand{\gx}[1]{\mbox{\textbf{#1}}}
\newcommand{\gxD}[1]{\mbox{\textbf{#1}}}
\newcommand{\gxS}[1]{\mbox{\textbf{\textsl{#1}}}}

\newcommand{\hilight}[1]{#1} 

\theoremstyle{plain}
\newtheorem{theorem}{Theorem}[section]
\newtheorem{proposition}[theorem]{Proposition}

\newtheorem{corollary}[theorem]{Corollary}
\newtheorem{conjecture}[theorem]{Conjecture}
\newtheorem{definition}[theorem]{Definition}
\newtheorem{example}[theorem]{Example}
\newtheorem{algorithm}[theorem]{Algorithm}  

\theoremstyle{remark}
\newtheorem{remark}[theorem]{Remark}

\numberwithin{equation}{section}
\newcommand{\ms}{\begin{math}}
\newcommand{\me}{\end{math}}

\newcommand{\bP}{\mathbb{Z}_+} 
\newcommand{\bQ}{\mathbb{Q}}

\newcommand{\xx}{\mathbf{x}}
\newcommand{\bh}{\mathbf{h}}
\newcommand{\po}{\overline{h}}

\newcommand{\cA}{\mathcal{A}}

\newcommand{\cC}{\mathcal{C}}

\newcommand{\cI}{\mathcal{I}}
\newcommand{\cJ}{\mathcal{J}}
\newcommand{\cL}{\mathcal{L}}
\newcommand{\nqs}{\mathcal{S}^*}
\newcommand{\qs}{\mathcal{S}}

\newcommand{\fldk}{\textbf{k}}

\newcommand{\Sym}{\textsl{Sym}}
\newcommand{\Qsym}{\textsl{QSym}}
\newcommand{\Nsym}{\textsl{NSym}}

\newcommand{\ncm}{\mathbf{m}}
\newcommand{\ncM}{\mathbf{M}}
\newcommand{\Tdot}{\dot T}
\newcommand{\SRS}{S^{RS}}
\newcommand{\QSRS}{\qs^{RS}}

\newcommand\strongof[1]{{#1}^+}
\newcommand\partitionof[1]{\widetilde{#1}}
\newcommand\reverse[1]{{#1}^*}

\newcommand{\refines}{\preccurlyeq}   

\newcommand{\suchthat}{:}
\newcommand{\revcon}{\Subset}    
\newcommand{\cskew}{\negthinspace\sslash\negthinspace}    
\newcommand{\cequiv}{\stackrel{c}{\sim}}    
\newcommand{\nearcat}{\odot}
\newcommand{\setnot}{\!\setminus\!}
\newcommand{\spam}{\text{span}}
\newcommand{\cshape}{C\text{-}shape}
\newcommand{\pshape}{P\text{-}shape}

\newlength\cellsize 
\savebox2{%
\begin{picture}(15,15)
\put(0,0){\line(1,0){15}}
\put(0,0){\line(0,1){15}}
\put(15,0){\line(0,1){15}}
\put(0,15){\line(1,0){15}}
\end{picture}}
\newcommand\cellify[1]{\def\thearg{#1}\def\nothing{}%
\ifx\thearg\nothing
\vrule width0pt height\cellsize depth0pt\else
\hbox to 0pt{\usebox2\hss}\fi%
\vbox to 15\unitlength{
\vss
\hbox to 15\unitlength{\hss$#1$\hss}
\vss}}
\newcommand\tableau[1]{\vtop{\let\\=\cr
\setlength\baselineskip{-16000pt}
\setlength\lineskiplimit{16000pt}
\setlength\lineskip{0pt}
\halign{&\cellify{##}\cr#1\crcr}}}
\savebox3{%
\begin{picture}(15,15)
\put(0,0){\line(1,0){15}}
\put(0,0){\line(0,1){15}}
\put(15,0){\line(0,1){15}}
\put(0,15){\line(1,0){15}}
\end{picture}}
\newcommand\expath[1]{%
\hbox to 0pt{\usebox3\hss}%
\vbox to 15\unitlength{
\vss
\hbox to 15\unitlength{\hss$#1$\hss}
\vss}}
\newcommand\bas[1]{\omit \vbox to \cellsize{ \vss \hbox to \cellsize{\hss$#1$\hss} \vss}}

\newcommand{\emt}{\mbox{ }}
\newcommand{\gst}{\mbox{\textbf{\textup{*}}}}
\newcommand{\gstD}{\gxD{*}}

\begin{document}

\bibliographystyle{hsiam}

\title[Skew quasisymmetric and noncommutative Schur functions]{Skew quasisymmetric Schur functions and noncommutative Schur functions}

\author{C. Bessenrodt}
\address{Institut f\"{u}r Algebra, Zahlentheorie und Diskrete Mathematik,
Leibniz Universit\"{a}t, 
D-30167 Hannover, Germany}
\email{\href{mailto:bessen@math.uni-hannover.de }{bessen@math.uni-hannover.de}}

\author{K. Luoto}
\address{Department of Mathematics, University of British Columbia, Vancouver, BC V6T 1Z2, Canada}
\email{\href{mailto:kwluoto@math.ubc.ca}{kwluoto@math.ubc.ca}}

\author{S. van Willigenburg}
\address{Department of Mathematics, University of British Columbia, Vancouver, BC V6T 1Z2, Canada}
\email{\href{mailto:steph@math.ubc.ca}{steph@math.ubc.ca}}
\thanks{
The second and third authors were supported in part by the National Sciences and Engineering Research Council of Canada. The third author was supported in part by the Alexander von Humboldt Foundation.}

\begin{abstract}
Recently a new basis for the Hopf algebra of quasisymmetric functions $QSym$, called quasisymmetric Schur functions, has been introduced by Haglund, Luoto, Mason, van Willigenburg. In this paper we extend the definition of quasisymmetric Schur functions to introduce skew quasisymmetric Schur functions. These functions include both classical skew Schur functions and quasisymmetric Schur functions as examples, and give rise to a new poset $\mathcal{L}_C$ that is analogous to Young's lattice. We also introduce a new basis for the Hopf algebra of noncommutative symmetric functions $NSym$. This basis of $NSym$ is dual to the basis of quasisymmetric Schur functions and its elements are the pre-image of the Schur functions under the forgetful map $\chi: NSym \rightarrow Sym$.  We prove that the multiplicative structure
constants of the noncommutative Schur functions, equivalently the
coefficients of the skew
quasisymmetric Schur functions when expanded in the quasisymmetric Schur basis, are nonnegative integers, satisfying a Littlewood-Richardson rule analogue that reduces to the classical Littlewood-Richardson rule under $\chi$.

As an application we show that the morphism 
\hilight{of algebras} 
from the algebra of Poirier-Reutenauer to $Sym$ factors through $NSym$. We also extend the definition of Schur functions in noncommuting variables of Rosas-Sagan in the algebra $NCSym$ to define quasisymmetric Schur functions in the algebra $NCQSym$. We prove these latter functions refine the former and their properties, and project onto quasisymmetric Schur functions under the forgetful map. Lastly, we show that by suitably labeling $\mathcal{L} _C$, skew quasisymmetric Schur functions arise in the theory of Pieri operators on posets.
\end{abstract}

\maketitle
\tableofcontents

\section{Introduction}\label{sec:intro}

At the beginning of the last century, Schur \cite{schur} identified   functions that would later bear his name as   characters of the irreducible
polynomial
representations of $GL(n, \mathbb{C})$. These functions   subsequently rose further in importance due to their ubiquitous nature. For example, in combinatorics they are the generating functions for semistandard Young tableaux, while in the representation theory of the symmetric group they form the image of the irreducible characters under the characteristic map. However, one of their most significant impacts has been as an orthonormal basis for the graded Hopf algebra of symmetric functions, $Sym$. More precisely, given partitions
$\lambda, \mu$,
the expansion of the product of Schur functions $s_\lambda , s_\mu$ in this basis is
$$ s_ \lambda s_\mu= \sum _ \nu c^ \nu _{\lambda\mu} s_ \nu, $$
where the $c^ \nu _{\lambda\mu}$ are known as Littlewood-Richardson coefficients.
However, this is not the only instance of Littlewood-Richardson coefficients.
In the
ordinary
representation theory of the symmetric group, taking the induced tensor product of   Specht modules $S^ \lambda $ and $S^\mu$ results in
$$\left(S^ \lambda \otimes S^\mu\right)\uparrow ^{S_n}=\bigoplus_\nu c^\nu _{\lambda\mu} S^ \nu .$$
Additionally, considering the cohomology $H^*(Gr(k,n))$ of the Grassmannian, the cup product of Schubert classes $\sigma _ \lambda$ and $\sigma_\mu$ satisfies
$$\sigma _ \lambda \cup \sigma _\mu = \sum_\nu c^ \nu _{\lambda\mu} \sigma _{\lambda}.$$

The $c^ \nu _{\lambda\mu}$ also arise in the expansion of skew Schur functions $s _{\nu/\mu}$ expressed in terms of Schur functions
$$s_ {\nu /\mu} = \sum_\lambda c^ \nu _{\lambda\mu} s_\lambda .$$
Skew Schur functions are themselves of importance, arising in discrete geometry as the weight enumerator of certain posets \cite{gessel},
 in the study of the general linear Lie algebra \cite{uglov-1},
and in mathematical physics in relation to spectral decompositions \cite{kirillov-1}.
Furthermore, Littlewood-Richardson coefficients also play an important role in several applications, such as in proving Horn's conjecture \cite{knutson-tao-1}.
More details on Littlewood-Richardson coefficients can be found in \cite{stanley-ec2}.

Therefore, the efficient computation of   Littlewood-Richardson coefficients is  a central problem, and to date their
computation falls mainly into two categories -- their precise computation, and relations they satisfy.
Regarding their computation, a combinatorial rule known as the Littlewood-Richardson rule (conjectured in \cite{LR-1} and proved in \cite{schutz-1,gthomas}) exists, and over the years a variety of reformulations have arisen in order to make their computation more straightforward  including \cite{berenstein-zel,fomin-greene-2,HLMvW-2}.
Meanwhile, regarding relations they satisfy, in  \cite{knutson-tao-1} it was shown that $c^{m \nu} _{m \lambda, m\mu} \geq  c^ \nu_{\lambda\mu}$, and further polynomiality properties were established in
 \cite{derksen-weyman-1,king-2}.
Furthermore, instances when they equate to 0 were identified in \cite{purbhoo-1} and when they equate to each other has been investigated in
\cite{billera-thomas-vW,gutschwager,mcnamara-vanw-1,RSvW, vW}.

As a consequence of the impact Schur functions have on other areas, and the combinatorial nature of the Littlewood-Richardson rule, Schur functions have been generalized to a number of analogues in the hope that these generalizations will also afford  combinatorial 
formulas to solve problems in related areas. Examples of analogues include Schur $P$ functions arising in the representation theory of the
double cover of the
symmetric group \cite{ macdonald-1, stembridge-1}, $k$-Schur functions connected to  the enumeration of Gromov-Witten invariants \cite{lapointe-morse}, cylindric Schur functions  \cite{mcnamara-1}, shifted Schur functions related to the representation theory of  ${\rm GL}(n)$ \cite{okounkov-1}, and factorial Schur functions that are special cases of double Schubert polynomials \cite{lascoux-1,molev-sagan}.
In addition, $Sym$ itself has been generalized: two of the most important generalizations being a nonsymmetric analogue and a noncommutative analogue, known as $QSym$ and $NSym$ respectively.

The nonsymmetric analogue $QSym$ is the Hopf algebra of quasisymmetric functions, and since its introduction as a source of generating functions for P-partitions \cite{gessel} quasisymmetric functions have been identified as generating functions for flags in  graded posets \cite{ehrenborg-1} and matroids \cite{billera-jia-reiner};
were shown to be the terminal object in the category of certain graded Hopf algebras \cite{aguiar-bergeron-sottile};
contain functions dual to the \textbf{cd}-index studied by discrete geometers \cite{billera-hsiao-vW};
arise as characters of a degenerate quantum group \cite{hivert-1}; investigate the behavior of random permutations \cite{stanley-2};
in addition to simplifying the calculation of symmetric functions such as Macdonald polynomials \cite{HHL-1, HLMvW-1} and Kazhdan-Lustig polynomials \cite{billera-brenti}.

Dual to $QSym$ is a noncommutative analogue $NSym$, the Hopf algebra of noncommutative symmetric functions first studied extensively in \cite{GKLLRT}.
In \cite{GKLLRT}, they proved $NSym$ is anti-isomorphic to Solomon's descent algebra \cite{solomon}, which in turn is anti-isomorphic to the dual of $QSym$  \cite{gessel, malvenuto-reut-1} and arises in the study of riffle shuffles \cite{bayer-diaconis,fulman}, and the study of Lie algebras \cite{garsia-reut,schocker-1}.
Meanwhile,
in representation theory, $NSym$ plays a role in the representation theory of the $0$-Hecke algebra, the $0$-quantum $GL_n$, and a $0$-quantized enveloping algebra \cite{ncsf4,ncsf5}.

Subsequently, these nonsymmetric and noncommutative analogues gave rise to nonsymmetric and noncommutative analogues of Schur functions. In $\Qsym$ the basis of fundamental quasisymmetric functions is often considered to form an analogue due to the aforementioned occurrence in the representation theory of a degenerative quantum group. Another analogue is  the basis of quasisymmetric Schur functions studied in \cite{HLMvW-1,HLMvW-2} that refine many classical combinatorial properties of Schur functions,
\hilight{although not yet the Littlewood-Richardson rule as the product of two quasisymmetric Schur functions often produces negative structure constants \cite[Section 7.1]{HLMvW-1}.
These functions have also}  
recently been applied in \cite{lauve-mason} to confirm a conjecture of Bergeron and Reutenauer.
In $\Nsym$ the basis dual to the fundamental quasisymmetric functions, known as the noncommutative ribbon Schur functions, are similarly considered to be Schur function analogues.
However, these are not the only noncommutative analogues that exist.
Other analogues include the noncommutative Schur functions of Fomin and Greene \cite{fomin-greene-1}, the Schur functions in noncommuting variables in $NCSym \subset NCQSym$ \cite{rosas-sagan}, and the free Schur functions arising in the algebra of Poirier and Reutenauer, $PR$ \cite{poirier-reutenauer}, sometimes called $FSym$ \cite{ncsf6}.
In this paper we propose a \emph{new} analogue for noncommutative Schur functions that differs from the previous analogues proposed,
viz. the basis of $\Nsym$ dual to the quasisymmetric Schur functions introduced in \cite{HLMvW-1}.
We establish connections between the various analogues in Section~\ref{sec:conc}.

For the moment, we give the connections between the algebras $Sym$, $QSym$, $NSym$, $NCSym$, $NCQSym$ and $PR$ below. Throughout we use $\chi$ to denote the forgetful map which allows algebra elements to commute.
$$
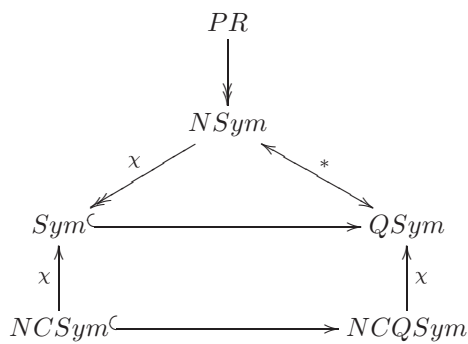
\xymatrix{&PR \ar@{->>}[d] &\\
&NSym \ar@{->>}[ld]_\chi \ar@{<->}[rd] ^\ast&\\
Sym   \ar@{<-}[d] _\chi  \ar@{^{(}->}[rr] && QSym \ar@{<-}[d] ^\chi \\
NCSym  \ar@{^{(}->}[rr]  && NCQSym}$$

This paper is structured as follows.
The intent of Section \ref{sec:back} is to provide sufficient background material to state and discuss our main results.
This includes defining new composition-indexed analogues of extant  partition-indexed combinatorial objects, such as tableaux, that arise in the literature of symmetric functions.
We also define skew quasisymmetric Schur functions and noncommutative Schur functions in Definition~\ref{def:skewqs}.   
In Section \ref{sec:NSymLR} we prove a combinatorial formula for skew quasisymmetric Schur functions in Proposition~\ref{prop:qschur-skew}; we state our main result, Theorem~\ref{thm:qschur-skew-pos}, a noncommutative Littlewood-Richardson rule, which gives a combinatorial interpretation to the multiplicative structure constants of the noncommutative Schur functions, thereby showing that these constants are nonnegative integers; and we discuss some  consequences of this rule such as recovering the classical Littlewood-Richardson rule, and noncommutative Pieri rules.
Section \ref{sec:NSymLRproof}  is devoted to the proof of the main result.

Section \ref{sec:apps}  provides some applications of quasisymmetric and noncommutative Schur functions, explicating the connections diagrammed above.
In Subsection~\ref{subsec:sym-SQS} we discuss a large class of skew quasisymmetric functions which are symmetric. This class includes the classical skew Schur functions.
In Subsection~\ref{subsec:PR} we demonstrate a new map showing that $\Nsym$, as an algebra, is a quotient of $PR$.
In Subsection~\ref{subsec:NCSym} we describe quasisymmetric Schur function analogues in $NCQSym$ which decompose the Schur functions analogues of \cite{rosas-sagan} in $NCSym$.
Lastly, in Subsection~\ref{subsec:pieri} we show how the skew quasisymmetric Schur functions can be interpreted from the
 viewpoint of Pieri operators  \cite{bergeron-vW} just as  the skew noncommutative Schur functions of Fomin and Greene can be.
Section  \ref{sec:conc}  briefly discusses future directions for research.

\medskip {\bf Acknowledgments.}
The authors would like to thank Sami Assaf, Sergey Fomin, Jim Haglund, Sarah Mason, Jean-Christophe Novelli, Mercedes Rosas, and Mike Zabrocki, for helpful discussions and suggestions that sparked fruitful paths of investigation.
\hilight{The authors would also like to thank Vasu Tewari and the referee for thoughtful comments.} 

\section{Background}\label{sec:back}

\subsection{Compositions, partitions, and tableaux}\label{subsec:compositions}
A \emph{weak composition} is a finite sequence of nonnegative integers, whose elements we call its \emph{parts}.
A \emph{strong composition}, or simply a \emph{composition}, is a finite sequence of positive  integers. (Thus every composition is a weak composition, but not vice versa.)
Given the weak composition $\alpha = (\alpha_1,\ldots, \alpha_k)$, we define its \emph{weight} as $|\alpha|= \alpha_1 +\cdots+ \alpha_k$ and its \emph{length} as $\ell(\alpha)=k$. If $|\alpha|=n$, we also write $\alpha\vDash n$.
There is a natural bijection between compositions $\alpha\vDash n$ and subsets of $[n-1]=\{1,\ldots,n-1\}$ which maps a composition to the set of its partial sums, not including $n$ itself, that is
\[  set(\alpha) := \{\alpha_1,\, \alpha_1 +\alpha_2,\,  \alpha_1 +\alpha_2+\alpha_3,\,  \ldots,\,n-\alpha_{\ell(\alpha)} \}. \]
\hilight {Following the convention of \cite{macdonald-1},} we 
say that $\beta$ \emph{refines} $\alpha$, denoted $\beta\refines\alpha$,
 if $|\alpha|=|\beta|$ and
\hilight{if we can obtain the parts of $\alpha$ by adding together consecutive parts of $\beta$.}
The \emph{reversal} of  $\alpha$, denoted $\reverse{\alpha}$, is the weak composition obtained by writing the parts of $\alpha$ in reverse order.
The \emph{underlying strong composition} of a weak composition $\gamma$, denoted $\strongof{\gamma}$, is the composition obtained by removing the zero-valued parts of $\gamma$ while keeping the nonzero parts in their same relative order.
A \emph{partition} is a composition whose parts are weakly decreasing.
If $\lambda$ is a partition with $|\lambda|=n$, we write $\lambda\vdash n$.
The \emph{underlying partition} of a weak composition $\alpha$, denoted $\partitionof{\alpha}$, is the partition obtained by sorting the nonzero parts of $\alpha$ into weakly decreasing order.
The \emph{empty composition (partition)}, denoted $\emptyset$, is the unique composition with weight and length zero.
The \emph{concatenation} of $\alpha=(\alpha_1,\ldots,\alpha_k)$ and $\beta=(\beta_1,\ldots,\beta_\ell)$ is $\alpha \beta =(\alpha_1,\ldots,\alpha_k,\beta_1,\ldots,\beta_\ell)$, while their \emph{near concatenation} is $\alpha\nearcat\beta =(\alpha_1,\ldots,\alpha_k+\beta_1,\ldots,\beta_\ell)$.
\begin{example}
For the weak composition, $\gamma=(2,3,0,1,4,0,2)$, we have $|\gamma|=12$, $\ell(\gamma)=7$, $\reverse{\gamma}=(2,0,4,1,0,3,2)$, $\strongof{\gamma}=(2,3,1,4,2)$.
For the composition $\alpha=(2,3,1,4,2)$, we have  $\partitionof{\alpha}=(4,3,2,2,1)$, and $set(\alpha)=\{2,5,6,10\}\subset [11]$.
Also, $(4,0,3,2,1,2)\refines(7,2,3)$,
 the concatenation $(1,2,3)(4,5)=(1,2,3,4,5)$, and $(1,2,3)\nearcat(4,5)=(1,2,7,5)$.
\end{example}

Given a composition $\alpha$, we say that the \emph{diagram of straight shape $\alpha$} is the left-justified arrangement of rows of \emph{cells}, where, following the English convention, the first (top) row of the diagram contains $\alpha_1$ cells, the second contains $\alpha_2$ cells, etc.
Viewing the diagram as a subset of $\bP\times\bP$,
we use $(row,column)$ pairs to index cells of a composition diagram, where row and column numbers start with 1.
We use the same symbol $\alpha$ to denote both a diagram and its shape when the usage is clear from context.
\begin{example}
\[ \tableau{\emt&\emt&\emt&\emt\\ \emt&\emt&\emt\\ \emt&\emt\\ \emt&\emt\\ \emt} \qquad \qquad
\tableau{\emt&\emt\\ \emt&\emt&\emt&\emt\\ \emt\\ \emt&\emt&\emt\\ \emt&\emt}
\]
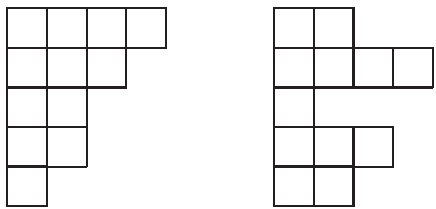
\centerline{Diagrams of the partition (4,3,2,2,1) and the composition (2,4,1,3,2)  }
\end{example}

\subsubsection{Poset of compositions} \label{subsec:comp-poset}

We say that the composition $\alpha$ is \emph{contained} in the composition $\beta$, denoted $\alpha \subset \beta$, if and only if $\ell(\alpha)\leq\ell(\beta)$ and $\alpha_i\leq \beta_i$ for all $1\leq i\leq \ell(\alpha)$.
We write $\alpha\revcon\beta$ to mean $\reverse{\alpha}\subset\reverse{\beta}$.

Young's lattice, which we denote $\cL_Y$, is the set of all partitions partially ordered by containment.
The empty partition is the unique minimal element  of $\cL_Y$.
In this lattice, $\nu$ covers $\mu$ if and only if $\nu $ can be obtained from $\mu$ by either appending a new part of size 1, or by incrementing some part of $\mu$ by 1, specifically the first part (i.e., leftmost part, or uppermost row, in a diagram) of a given size.
We define an analogous partial order on the set of compositions, whose importance will be apparent in the sections that follow.

\begin{definition}[Composition poset] \label{def:comp-poset}
We say that
the
composition  \emph{$\gamma$ covers $\beta$}, denoted $\beta\lessdot_C \gamma $, if $\gamma $ can be obtained from $\beta $ either by prepending $\beta $ with a new part of size 1, or by adding 1 to the first (leftmost) part of $\beta$ of size $k$ for some $k$.
The partial order $\leq_C$ defined on the set of all compositions is the transitive closure of these cover relations, and the resulting poset we denote $\cL_C$.
\end{definition}

\begin{remark} \label{rem:LC}
Note that $\beta\lessdot_C \gamma $ implies $\beta\revcon\gamma$, but not vice versa.
Clearly $\cL_C$ is graded by $rank(\alpha)=|\alpha|$ and has a unique minimal element $\emptyset$.
However $\cL_C$ is not a lattice; neither meets nor joins are defined in general.
\end{remark}

\begin{example}
The composition $(2,3,2)$ is covered by $(1,2,3,2)$, $(3,3,2)$ and $(2,4,2)$, but not by $(2,3,3)$.
The compositions $(2,2,2)$ and $(2,3,2)$ do not have a meet
as they both lie over $(1,2,2)$ as well as $(1,2,1)$.
The compositions $(2,3,1)$ and $(2,3,4)$ do not have a join as they both lie below $(4,4,4)$ as well as $(4,5,4)$.
\end{example}

\subsubsection{Skew shapes and tableaux}  \label{sec:skew-shapes} \label{sec:Tab-defs}

We extend the notions of shapes and diagrams to the skew case.
A diagram of skew shape is indexed by an ordered pair of compositions, but in contrast to diagrams of straight shape, we must distinguish between \emph{skew partition shapes $\nu/\mu$}, and \emph{skew composition shapes $\gamma\cskew\beta$}, as defined below.
In both cases, as with diagrams of straight shape, we use the same symbol to denote both a diagram (a configuration of cells) and its shape (an ordered pair of compositions) when the meaning is clear.

\begin{definition}[Skew shapes]
Given partitions $\mu \subset \nu$, the \emph{diagram of skew partition shape $\nu/\mu$} comprises  those cells in the diagram of shape $\nu$ that are not in the diagram of shape $\mu$ when the diagram of $\mu$ is positioned in the upper left of that of $\nu$.
\hilight{We write $|\nu/\mu|:= |\nu|-|\mu|$.} 

Given compositions $\beta \revcon \gamma$,
the \emph{diagram of skew composition shape $\gamma \cskew \beta$} comprises  those cells in the diagram of shape $\gamma $ that are not in the diagram of shape $\beta $ when the diagram of $\beta$ is positioned in the lower left of that of $\gamma$.
\hilight{We write  $|\gamma\cskew\beta|:= |\gamma|-|\beta|$.} 
\end{definition}
\begin{example}
\[  \tableau{\gst& \gst&\emt&\emt\\ \gst& \gst&\emt\\ \gst&\emt\\ \gst&\emt\\ \emt}
\qquad \qquad \tableau{\emt&\emt\\ \gst&\emt&\emt&\emt\\ \gst\\ \gst&\gst&\emt\\ \gst&\gst} \]
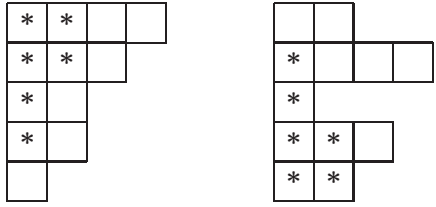
\centerline{Diagrams of $(4,3,2,2,1)/(2,2,1,1)$  and $(2,4,1,3,2)\cskew(1,1,2,2)$}
\[ |(4,3,2,2,1)/(2,2,1,1)| \;=\; |(2,4,1,3,2)\cskew(1,1,2,2)| \;=\; 6 \]
\end{example}

Note that under this definition, a straight shape is a skew shape of the form $\lambda/\emptyset$ or $\alpha\cskew\emptyset$ respectively.

\begin{definition}[Strips]
A \emph{vertical strip} is a skew shape (either partition or composition) whose diagram contains at most one cell per row.
A \emph{horizontal strip} is a skew shape whose diagram contains at most one cell per column.
\hilight{}
\end{definition}

A \emph{partition-shaped tableau}  is a filling $T\colon \nu/\mu\to\bP$ of the cells of a (skew) partition diagram with positive integers.
A \emph{semistandard Young tableau}  (SSYT) is a partition-shaped tableau in which
the entries in each row are weakly increasing from left to right, and
the entries in each column are strictly increasing from top to bottom.
A \emph{standard} Young tableau (SYT) is an SSYT in which
the filling is a bijection $T\colon \nu/\mu\to[n]$ where $[n]=\{1,2,\ldots,n\}$ and $n=|\nu/\mu|$.

A \emph{(semi-)standard reverse tableau} (SSRT or SRT) is like an SSYT or SYT except that we reverse the inequalities:
the entries in each row are weakly decreasing from left to right, and
the entries in each column are strictly decreasing from top to bottom.
All concepts relating to Young tableaux have their reverse tableau counterparts.
In this article we primarily make use of reverse tableaux for partition shapes,
for consistency with \cite{HLMvW-1,HLMvW-2}, 
\hilight{and moreover to simplify the proofs by avoiding the introduction of additional notation and machinery.}

\begin{definition}[Composition tableau] \label{def:SSCT}
Given compositions $\beta\revcon\gamma$,  consider the cells of their respective diagrams as subsets of $\bP\times\bP$ indexed by $(row,column)$, arranged according to the convention for the skew diagram of shape $\gamma\cskew\beta$, so that the last row of the cells of $\beta$ lie in the last row $\ell(\gamma)$ of cells of $\gamma$.

We say that the cell $(i,k)$ \emph{attacks} the cell $(j,k+1)$ if $i<j$, $(j,k+1)\in \gamma\cskew\beta$, and $(i,k+1)\notin \beta$, although possibly $(i,k)\in \beta$.

A filling $T: \gamma\cskew \beta\to \bP$  is a \emph{semistandard composition tableau $T$  (SSCT) of shape $\gamma\cskew\beta$} if it satisfies the following conditions:
\begin{enumerate}
\item
Row
entries are weakly decreasing from left to right.
\item The entries in the first column are strictly increasing from top to bottom.
\item  If
$(i,k)\in \gamma$ attacks  $(j,k+1)$ and either  $(i,k)\in \beta$ or $T(j,k+1)\leq T(i,k)$, then $(i,k+1)\in \gamma\cskew \beta$ and  $T(j,k+1)< T(i,k+1)$.
\end{enumerate}
We say that $T$ is  \emph{standard} (an SCT) if $T$ is injective and its range is precisely $[n]$, where $n= | \gamma\cskew\beta|$.
\end{definition}
\begin{example}
\[ \begin{array}{ccc}
  \tableau{ \gst & \gst & \gst &2 \\ \gst & \gst &6 \\ \gst &9 &4 \\ 8 &5 &1\\ 6 &4 \\ 3 &3 }
  &\text{ }\qquad &
  \tableau{ 3 &3 &1\\  6 &5&4 &2\\ 8 &4 \\  \gst & \gst &6 \\ \gst & \gst & \gst \\ \gst &9  }
  \\  \\ \text{an } SSRT \text{ of shape }  & & \text{an } SSCT \text{ of shape }
  \\   (4,3,3,3,2,2)/(3,2,1) & & (3,4,2,3,3,2)\cskew(2,3,1)
\end{array} \]
\end{example}

This definition of SSCT is consistent with that of \cite{HLMvW-1} for straight composition shapes.
We let $SSRT(\nu/\mu)$ (resp. $SRT(\nu/\mu)$) denote the set of all SSRT (resp. SRT) of shape~$\nu/\mu$.
Similarly, we let $SSCT(\gamma\cskew\beta)$ (resp. $SCT(\gamma\cskew\beta)$) denote the set of all SSCT (resp. SCT) of shape $\gamma\cskew\beta$.
We write $sh(T)$ to denote the \emph{shape} of a tableau:  $sh(T)=\nu/\mu$ if $T\in SSRT(\nu/\mu)$, or  $sh(T)= \gamma\cskew\beta$ if $T\in SSCT(\gamma\cskew\beta)$.

Recall that a \emph{saturated chain} in a poset is a (finite) sequence of 
\hilight{consecutive} 
cover relations.
There is a well-known natural bijection between SYT (equivalently, SRT) and saturated chains in Young's lattice $\cL_Y$.
Likewise we have the following.

\begin{proposition} \label{prop:sct-chains}
There is a natural bijection between $SCT(\gamma\cskew\beta)$ and the set of saturated chains in $\cL_C$ from $\beta$ to $\gamma$.
\end{proposition}
\begin{proof}
Suppose that  $|\gamma\cskew\beta|=n$ and $T\in SCT(\gamma\cskew\beta)$.
Define the sequence of compositions
\begin{equation}  \label{eqn:C-chain}
 \alpha^n=\beta\revcon\alpha^{n-1}\revcon\cdots\revcon\alpha^1\revcon\alpha^0=\gamma
\end{equation}
by the rule
\begin{equation} \label{eqn:T-chain}
\alpha^{k-1}=\alpha^k\cup T^{-1}(k) \qquad\text{ for all } 1\leq k\leq n.
\end{equation}
That is, $\alpha^{k-1}$ is obtained from $\alpha^k$ by adding to its diagram the cell position of $T$ that contains~$k$.  For example,
\[
T=\tableau{ 4&2\\ \gst &1 \\ \gst &\gst&\gst \\ \gst &\gst &3} \quad\longleftrightarrow\quad
(1,3,2)\revcon(1,1,3,2)\revcon(1,1,3,3)\revcon(2,1,3,3)\revcon(2,2,3,3).
\]
We claim that this rule defines the desired bijection.
To see that the sequence is a chain in~$\cL_C$, proceed by induction on $k$, starting with $k=n$.
Since the entries in the first column of $T$ are increasing top to bottom,
if the entry $k$ appears in the first column of $T$, then all cells below it in the same column of $T$ already belong to $\alpha^k$, hence $\alpha^{k-1}$ is obtained by prepending a new part of size 1, and so $\alpha^k\lessdot_C\alpha^{k-1}$.
Otherwise, since row entries are decreasing, the cell $T^{-1}(k)=(i,j+1)$ appears immediately to the right of either a higher numbered cell or a cell in $\beta$, either of which by hypothesis belongs to $\alpha^k$, so $\alpha^{k-1}$ is obtained by appending a new cell to the end of some row of $\alpha^k$
(i.e., incrementing some part of $\alpha^k$).
Moreover, there can be no higher row $i'<i$ of $\alpha^k$ of length $j$ 
since otherwise in $T$,  cell $(i',j)$ would attack $(i,j+1)$ with either $(i',j)\in\beta$ or  $T(i,j+1)<T(i',j)$, and with either $(i',j+1)\notin\gamma$ or $T(i,j+1)>T(i',j+1)$, contradicting that $T$ is an SCT.
So $i$ is the highest row of $\alpha_k$ having length $j$, and again $\alpha^k\lessdot_C\alpha^{k-1}$ as desired.
Thus every SCT determines a unique saturated chain.

Conversely, suppose we have a saturated chain in $\cL_C$ of the form \eqref{eqn:C-chain}.
Let $T$ be the filling of $\gamma\cskew\beta$ determined by the relations  \eqref{eqn:T-chain}.
Then $T:\gamma\cskew\beta\to[n]$ is a bijection.
The covering relations ensure that the first column of $T$ is increasing and that all rows are decreasing.
Suppose cell $(i',j)$ attacks $(i,j+1)$ in $T$ with $T(i,j+1)=k$ and either $(i',j)\in\beta$ or, say, $T(i,j+1)<T(i',j)=k'$.
Then in either case we have $(i',j)\in\alpha^k$.
Since $\alpha^{k} \lessdot_C \alpha^{k-1}$, row $i$ is the highest row of $\alpha_k$ of length $j$, hence $(i',j+1)\in\alpha^k$ and so $(i',j+1)=T^{-1}(k'')$ for some $k<k''$, that is, $T(i,j+1)=k<k'' =T(i',j+1)$.
Since we considered arbitrary attacking cells, $T$ satisfies the conditions of an SCT.
\end{proof}

After reading the above proof, the equivalent bijection between $SRT(\nu/\mu)$ and the set of saturated chains in $\cL_Y$ from $\mu$ to $\nu$ should be clear.

\subsubsection{Tableau properties} \label{sec:Tab-props}

Unless otherwise indicated, the definitions in this section apply to both reverse tableaux (SSRT) and composition tableaux (SSCT), and to those of skew shape as well as straight shape.

\smallskip
The \emph{content} of a tableau $T$, denoted $cont(T)$, is the weak composition $\tau$  where $\tau_i$ denotes the number of entries of $T$ with value $i$.
The \emph{column word of $T$}, denoted $w_{col}(T)$, is the word consisting of the entries of each column of $T$ arranged in increasing order, beginning with the first (leftmost) column.

\begin{example}
\[ T = \tableau{ 4 &3&2&1 \\ 6&2&1 \\ \gst &1 \\ \gst&\gst&4&4&2 }
\qquad
\begin{array}{rcl} \\
cont(T)       &=& (3,3,1,3,0,1) \\ \\
 w_{col}(T) &=& 46\,123\,124\,14\,2 \\
\end{array}
\]
\end{example}

\textbf{Note: } The column word defined here should \emph{not}
be confused with the column reading word used in the papers \cite{HLMvW-1,HLMvW-2,mason-1}.
However, this definition of column word is consistent with the usual definition of column reading word for SSRT.

Let $T$ be a standard tableau, containing $n$ cells.
The \emph{descent set} of $T$, denoted $descents(T)$, is the set of those entries $i\in[n-1]$ such that $i+1$ is in a column weakly to the right of $i$ in $T$.
The \emph{descent composition} of $T$, denoted $Des(T)$, is the composition of $n$ associated to $descents(T)$ via partial sums, that is, $set(Des(T))=descents(T)$.

Given a composition $\alpha$, the \emph{canonical composition tableau of shape $\alpha$}, denoted  $U_\alpha$ is the unique SCT of shape $\alpha$ whose descent composition is $\alpha$, that is, $Des(U_\alpha)=\alpha$.
One can construct $U_\alpha$ by starting with the unfilled diagram of shape $\alpha$ and  consecutively numbering the cells  in the last (bottom) row  from left to right, then the next to last row, etc. in decreasing fashion.
\begin{example}
\[ T=\;\; \tableau{ 7  &4&2&1 \\  8 & 3\\ \gst \\\gst&\gst&6&5}\quad \quad
U_{1314} = \;\; \tableau{ 1 \\ 4 & 3 & 2 \\ 5 \\ 9 & 8 & 7 & 6 }
\qquad
\begin{array}{l} \\
descents(T)=\{3,4,7\}
 \\ Des(T)=(3,1,3,1)
 \\ \\ Des(U_{1314})=(1,3,1,4)
\end{array}
\]
\end{example}

\medskip
Let $T$ be a filling of a skew shape (composition or partition) such that the entries in each column are distinct.
(Such fillings include SSRT, and also SSCT as is straightforward to verify from property (3) of Definition \ref{def:SSCT}.)
The \emph{standard order $<_T$ of cells determined by $T$}  is the total order on the cells of the skew shape given by
\[ (i,j) <_T (i',j') \quad \text{ if }\;  T(i,j) < T(i',j')\; \text{ or } \left(\; T(i,j) = T(i',j') \;\text{ and }\; j>j' \;\right). \]
The \emph{standardization of $T$}, denoted $std(T)$, is the filling of $sh(T)$ obtained from $T$ by renumbering the cells in their standard order in consecutive increasing fashion.

The \emph{column sequence} of a standard tableau $T$, denoted $colseq(T)$,
is the sequence $(j_n, j_{n-1},$ $\ldots,  j_1)$ of the column numbers of the cells of $T$ listed in decreasing order of their entries, that is, $T^{-1}(k)=(i_k,j_k)$ for some row $i_k$ for all $1\leq k\leq n$.
\begin{example}
\[ T = \tableau{ 4 &2&1&1 \\ \gst &1\\\gst&\gst&4&2}
\qquad std(T)=  \tableau{ 7 &5&2&1 \\ \gst &3\\\gst&\gst&6&4}
\qquad colseq(std(T))=(1,3,2,4,2,3,4)
 \]
\end{example}
\medskip

The following  fact is known for SSYT, and the proof is equally straightforward for SSRT and SSCT.

\begin{proposition} \label{prop:standardization}
The standardization $std(T)$ of a tableau $T$ is a standard tableau.
There is a natural bijection between $SSCT(\gamma\cskew\beta)$  (resp. $SSRT(\nu/\mu)$), and ordered pairs $(\hat{T},\tau)$ where $\hat{T}\in SCT(\gamma\cskew\beta)$ (resp. $\hat{T}\in SRT(\nu/\mu)$) and $\tau$ is a weak composition such that $\tau\refines Des(\hat{T})$.
More precisely, $T\leftrightarrow (std(T),cont(T))$.
\end{proposition}
\begin{proof}
Let $T$ be a filling of  $\gamma\cskew\beta \vDash n$ such that the entries in each column are distinct, without assuming \emph{a priori} that $T$ is an SSCT, and let $\hat{T}=std(T)$.
Note that standardization preserves the standard order of cells in the skew shape.
With respect to this same standard order of cells, $T \colon \gamma\cskew\beta\to\bP$ is a weakly increasing function, and $\hat{T} \colon \gamma\cskew\beta\to[n]$ is strictly increasing.
Thus the entries of $T$ weakly decrease along a row left to right (i.e., the cells decrease in standard order) if and only if the entries of $\hat{T}$ do.
Likewise the entries of $T$ increase along the first column top to bottom if and only if the entries of $\hat{T}$ do.
Similarly, there are attacking cells $(i,j)$ and $(i',j+1)$ in $\hat{T}$ which violate condition $(3)$ of Definition \ref{def:SSCT} if and only if  the same cells of $T$  violate the condition.
Thus $\hat{T}$ is an SCT if and only if  $T$ is an SSCT.

Suppose $k\in descents(\hat{T})$. Let $\hat{T}^{-1}(k)=(i,j)$ and $\hat{T}^{-1}(k+1)=(i',j')$.
Now $k\in descents(\hat{T})$ implies $j\leq j'$, and $(i,j)<_{\widehat{T}} (i',j')$ implies $(i,j)<_{T} (i',j')$. These in turn imply $T(i,j)<T(i',j')$. Thus $cont(T)\refines Des(\hat{T})$.

Conversely, let $\tau\refines Des(\hat{T})$ and let $T'\colon \gamma\cskew\beta\to\bP$ be the uniquely determined filling which is weakly increasing with respect to the standard order of the cells determined by $\hat{T}$ and such that $cont(T')= \tau$.
This is always possible since $\tau\refines Des(\hat{T})$.
Then the entries in each column of $T'$ must be distinct, $\hat{T}=std(T')$, and by the above, $T'$ is an SSCT if and only if $\hat{T}$ is an SCT.

The case for showing that  $\hat{T}=std(T)$ is a SRT when $T$ is a $SSRT$, and the analogous bijection, is proved in
a similar fashion.
\end{proof}

\subsubsection{Mason's bijection $\rho$} \label{sec:mason}

Mason \cite{mason-1} described a natural bijection between reverse tableaux of straight shape and objects called \emph{semistandard skyline fillings}.
In subsequent work, the authors of \cite{HLMvW-1} extended this to a bijection $\rho:SSCT\to SSRT$ between composition tableaux and reverse tableaux of straight shape.
One characterization of the bijection is that $\rho$ preserves the underlying \emph{column tabloid}, i.e., the set of entries within each column.
Thus one direction is easy to compute:
 the reverse tableau $\rho(T)$ is obtained from the composition tableau $T$ by simply sorting the entries of each column and top-justifying the entries.
The opposite direction is only slightly harder.
Given a reverse tableau $T'$, compute the inverse image $T=\rho^{-1}(T')$ as follows.
Take the set of entries in the first column of $T'$ and write them in increasing order in the first column of $T$.
Then processing the remaining columns of $T'$ in left to right order, and the entries within each column in descending order,  each entry is placed as high as possible so as to still maintain weakly decreasing rows.

\begin{example}
\[ \tableau{8&6&5&4&3\\ 7&5&4&2\\ 4&3&2 \\ 2&2}
 \qquad \stackrel{\rho^{-1}}{\longrightarrow} \qquad
\tableau{2&2&2&2 \\ 4&3 \\ 7&6&5&4&3 \\ 8&5&4} \]
\end{example}

If  $T\in SSCT(\alpha)$,  then $\rho(T)\in SSRT(\partitionof{\alpha})$.
In general, any property of SSRT of straight shape that depends only on the column tabloid can be naturally extended to SSCT of straight shape via the bijection $\rho$.
The bijection $\rho$  extends to skew tableaux in the following sense.

\begin{proposition} \label{prop:rho-general}
Let $\beta $ be a composition, let  $\mu=\partitionof{\beta}$,
and let $\nu$ be a partition with $\mu\subset \nu $.
Let \[ \cC_{\nu,\beta} := \{ T \colon T \in SSCT(\gamma\cskew\beta) \text{ for some } \gamma >_C \beta \text{ with } \partitionof{\gamma}=\nu \}. \]
Then there is a natural bijection between $\cC_{\nu,\beta}$ and $SSRT( \nu/\mu )$ that preserves the set of entries in each column of paired tableaux.
This specializes to Mason's bijection $\rho$ in the case $\beta=\emptyset$.
\end{proposition}
\begin{proof}
In view of Proposition  \ref{prop:standardization}, it suffices to show the bijection between the subset $\cC'_{\nu,\beta}$ of  $\cC_{\nu,\beta}$ comprising the SCT of appropriate shapes, and $SRT(\nu/\mu)$.

The column sequence of a standard tableau determines the set of entries in each column, i.e., the tabloid. Conversely, the column sequence can be recovered from the tabloid.
An 
SRT  of shape $\nu/\mu$ can be recovered from its base shape $\mu$ and its column sequence since the column sequence specifies a unique maximal chain in the interval of $\cL_Y$ from $\mu$ to $\nu $.
Similarly, an 
SCT of shape $\gamma\cskew\beta$ can be recovered from its base shape $\beta$ and its column sequence since the column sequence specifies a unique maximal chain in the interval of $\cL_C$ from $\beta$ to $\gamma$.
From the cover relations of $\cL_Y$ and $\cL_C$, it is clear that a column sequence can be ``applied'' to a base composition $\beta$ to determine a saturated chain in $\cL_C$ ending in 
$\gamma$ 
if and only if the column sequence can be ``applied'' to the base partition $\partitionof{\beta}$ to determine a saturated chain in $\cL_Y$ ending in $\nu=\partitionof{\gamma}$.
Since the column sequence determines a unique chain in each of the respective posets, the map is bijective.
\end{proof}

\subsection{RSK correspondence} \label{sec:RSK}

A \emph{word} over an ordered \emph{alphabet} $A$ is a finite sequence of elements of $A$, the elements referred to as \emph{letters}.
For our purposes, $A$ is the set of positive integers $\bP$ with its natural ordering.
Instead of using the usual \emph{Robinson-Schensted-Knuth (RSK) correspondence} as described for example in \cite{fulton-1} or \cite{stanley-ec2}, here we will
use a ``reverse'' variant giving
a bijection between words $w$ and ordered pairs $(P(w),Q(w))$ of  SSRT of the same shape, where $Q(w)$ is standard;
this version is discussed in \cite{HLMvW-1}.
Recall that for the RSK correspondence, \emph{Schensted insertion} \cite{fulton-1} is used, i.e., 
an algorithm for inserting an integer $k$ into an SSRT  $S$ to obtain a new SSRT, denoted $S\gets k$, having one more cell than $S$, including the inserted integer as an entry.
We emphasize that in our context reverse bumping rules are used rather 
than doing ordinary Schensted insertion and then
reversing the tableaux.
Two words $w$ and $w'$ are said to be \emph{Knuth equivalent} or \emph{$P$-equivalent}, denoted $w\stackrel{P}{\sim}w'$, if $P(w)=P(w')$.
It is a fact that $P(w_{col}(T))=T$, that is, the column reading word of a straight SSRT belongs to its equivalence class.
The words $w$ and $w'$ are said to be  \emph{$Q$-equivalent}, denoted $w\stackrel{Q}{\sim}w'$, if $Q(w)=Q(w')$.
If $w$ and $w'$ are permutations, written in one-line notation and viewed as words,
then we also say that $w$ and $w'$ are  \emph{(resp. dual) Knuth equivalent} if they are (resp. $Q$-) $P$-equivalent.
Note that equivalence classes with respect to RSK and reverse RSK are the same.
These notions are discussed further in Section  \ref{sec:NSymLRproof}.

We define the \emph{rectification} of a skew (composition or reverse) tableau $T$ to be the straight composition tableau $rect(T):=\rho^{-1}(P(w_{col}(T)))$.
When we want to consider the rectification as a reverse tableau, we write $\rho(rect(T))$.
The \emph{rectified composition shape} of a skew 
tableau~$T$, denoted  $\cshape(T)$, is the shape of $rect(T)$.
The \emph{rectified partition shape} 
of~$T$, denoted  $\pshape(T)$, is the shape of $\rho(rect(T))$.
As  discussed below, it can be shown that
rectification preserves descents of tableaux, i.e., $Des(T) =Des(rect(T))$.
\begin{example}
\[  T=\; \tableau{ 4 &3&1\\ 8&6  \\ \gst &\gst &7 &5&2  \\ \gst &\gst &\gst \\ \gst &9 }
\quad
rect(T)=\; \tableau{ 4 & 3 &1\\ 8 & 7 &5 &2 \\ 9 & 6 }
\quad
\begin{array}{c}
\cshape(T) = (3,4,2) \\  Des(T) = Des(rect(T)) \;=\; (1,3,2,2,1)
\end{array}
 \]
\end{example}

\subsection{Symmetric and quasisymmetric functions}\label{subsec:sym-qsym}

The \emph{algebra of quasisymmetric functions}, $\Qsym$, is a subalgebra of
$\bQ[[X]]$, the formal power series ring over the commuting variables $X=\{x_1,x_2,\ldots\}$ indexed by the positive integers, graded by total monomial degree.
$\Qsym$ in turn contains the \emph{algebra of symmetric functions}, $\Sym$, as a subalgebra.
We direct the reader to \cite{stanley-ec2} for a basic introduction to the algebras $\Sym$ and $\Qsym$.

The basis elements of $\Sym$ are naturally indexed by partitions,
while the basis elements of $\Qsym$ are naturally indexed by compositions.
Both algebras have a natural monomial basis, which we denote $\{m_\lambda\}$ and $\{M_\alpha\}$ respectively.
$\Qsym$ has a second important basis known as the \emph{fundamental basis}, denoted $\{L_\alpha\}$.
One of the most important bases of $\Sym$ are the \emph{Schur functions}, denoted $\{s_\lambda\}$.
In \cite{HLMvW-1} the authors defined another basis of $\Qsym$, the \emph{quasisymmetric Schur functions}, denoted $\{\qs_\alpha\}$, which naturally refine the Schur functions.
We present below various known formulas for these bases.
First some notation.
Let $\xx^\gamma$ be the monomial indexed by the weak composition $\gamma$,
for example $\xx^{(3,0,1,2)} = x_1^3x_3x_4^2$.

\begin{alignat}{3}
M_\alpha &= \;\;\sum_{\strongof{\gamma}= \alpha} \xx^\gamma \\
 L_\alpha &= \;\;\sum_{\gamma\refines \alpha} \xx^\gamma &&= \;\;\sum_{\beta\refines\alpha} M_\beta   \\
\qs_\alpha &= \sum_{T\in SSCT(\alpha)} \xx^{cont(T)}
\;\; &&=  \sum_{T\in SCT(\alpha)} L_{Des(T)} \label{eqn:qschur} \\
m_\lambda &= \;\;\;\sum_{\partitionof{\gamma}=\lambda} \;\xx^\gamma
\;\; &&= \quad \;\sum_{\partitionof{\alpha}=\lambda} M_\alpha \label{eqn:mono}  \\
s_{\lambda} &= \sum_{T\in SSRT(\lambda)} \xx^{cont(T)}
\;\; &&= \sum_{T\in SRT(\lambda)} L_{Des(T)}
\;\; &&= \;\; \sum_{\partitionof{\alpha}=\lambda} \qs_\alpha \label{eqn:schur}
\end{alignat}

\subsection{Skew quasisymmetric Schur functions and noncommutative Schur functions}
\label{subsec:skewqs} \label{subsec:hopf}

The reader may find introductory material regarding Hopf algebras in \cite{hopf-intro, montgomery-1}.
An algebra $A$ over a field $\fldk$ can be thought of as a vector space with an associative product (bilinear map) $\,\cdot:A\otimes A \to A$ and a \emph{unit} (linear map) $u:\fldk\to A$.
A \emph{coalgebra} $C$ can be thought of as a vector space with a coassociative \emph{coproduct}  $\,\Delta:C\to C\otimes C$ and a \emph{counit} $\varepsilon:C\to\fldk$ (both linear maps).
A \emph{bialgebra} $H$ has both algebra and coalgebra structures satisfying certain compatibility conditions, such as $\Delta(a\cdot b) = (\Delta a)\cdot(\Delta b)$.
A \emph{Hopf algebra} is a bialgebra that has an automorphism known as an \emph{antipode}, which we will not describe here.
\hilight{Every connected, graded bialgebra is a Hopf algebra, for which we refer the reader to any of 
\cite[Section 2.3.3]{aguiar-mahajan-2}, \cite[Lemma 2.1]{ehrenborg-1}, or \cite[Proposition 8.2]{milnor-moore} for further details.}

\hilight{If $H=\bigoplus_{n\geq 0}H_n$ is a graded Hopf algebra of finite type, then $H^*=\bigoplus_{n\geq 0}H_n^*$ is the graded Hopf dual, where $H_n^*$ is the dual of $H_n$ respectively for each $n$.
Accordingly, there is a nondegenerate bilinear form $\langle\cdot,\cdot\rangle:H\otimes H^*\to \fldk$ that pairs 
the elements of any basis $\{B_i\}_{i\in\cI}$ of $H_n$ (for some index set $\cI$) and its dual basis $\{D_i\}_{i\in\cI}$ of $H_n^*$,
given by
$\langle B_i, D_j\rangle = \delta_{ij}$.
Duality is exhibited in that} the product coefficients of one basis are the coproduct coefficients of its dual basis and vice versa, i.e., \begin{eqnarray*}
B_i\cdot B_j = \sum_h a^h_{i,j} B_h &\qquad \Longleftrightarrow\qquad &
\Delta D_h = \sum_{i,j} a^h_{i,j} D_i \otimes D_j,  \\
D_i \cdot D_j = \sum_h b^h_{i,j} D_h &\qquad \Longleftrightarrow\qquad &
\Delta B_h = \sum_{i,j} b^h_{i,j} B_i \otimes B_j.
\end{eqnarray*}
As noted in Lam et al. \cite{LLS-hopfLR}, the coproduct allows one to define \emph{skew}
elements, indexed by ordered pairs of indices (usually written $i/j$), by
\begin{equation}  \label{eqn:generic-skew}
\Delta B_i = \sum_{j}  B_{i/j} \otimes B_j.
\end{equation}

$\Qsym$ has a Hopf algebra structure where the coproduct \cite{GKLLRT,malvenuto-reut-1} is given by
\begin{equation} \label{eqn:qscoprod}
 \Delta M_\alpha  = \sum_{\beta\gamma=\alpha} M_\beta\otimes M_\gamma,
\quad\text{equivalently,}\quad
\Delta L_\alpha  = \sum_{\beta\gamma=\alpha,\text{ or}\atop \beta\nearcat\gamma=\alpha} L_\beta\otimes L_\gamma.
\end{equation}
$\Sym$ has a Hopf algebra structure inherited from $\Qsym$.
The graded Hopf dual of $\Qsym$ is isomorphic to $\Nsym$, the \emph{algebra of noncommutative symmetric functions}  \cite{GKLLRT}, while $\Sym$ is self-dual as a Hopf algebra.
The duality pairing for $\Sym$ coincides with the (standard) Hall inner product, so that bases of $\Sym$ which are dual in the classical sense are also dual in the Hopf algebra sense.
Under this pairing, the monomial basis $\{m_\lambda\}$ is dual to the basis of complete symmetric functions $\{h_\lambda\}$, while the basis of Schur functions $\{s_\lambda\}$ is self-dual; specifically, each Schur function is dual to itself.
Equation \eqref{eqn:generic-skew} specializes to a formula for the classical \emph{skew Schur functions} in terms of the coproduct on $\Sym$ \cite{LLS-hopfLR}:
\begin{equation}
\label{eqn:skewSchurs}  \Delta s_\nu = \sum_{\mu}  s_{\nu/\mu} \otimes s_\mu.
\end{equation}
We use Equation \eqref{eqn:generic-skew} to define the \emph{skew quasisymmetric Schur functions}.
\begin{definition}[Skew quasisymmetric Schur functions]\label{def:skewqs}
The skew quasisymmetric Schur functions $\qs_{\gamma\cskew \beta}$ are defined implicitly
by the equations
\begin{equation}
\label{eqn:skewQSchurs}  \Delta \qs_\gamma = \sum_{\beta}  \qs_{\gamma\cskew \beta} \otimes \qs_\beta
\end{equation}
where $\beta$ ranges over all compositions.
\end{definition}

\medskip
A morphism of Hopf algebras $\phi:A\to B$ induces a morphism $\phi^*:B^*\to A^*$.
Suppose that $A$ and $B$ have bases $\{a_i\}_{i\in\cI}$ and $\{b_j\}_{j\in\cJ}$ respectively,
and that $\cJ$ can be partitioned into blocks (equivalence classes) $\cJ=\bigcup_{i\in\cI}\cJ_i$ such that $\phi(a_i)=\sum_{j\in\cJ_i} b_j$.
Then it can be shown that $\phi^*(b^*_j)=a^*_i$ for all $j\in\cJ_i$.

The inclusion $\Sym\hookrightarrow\Qsym$ therefore induces a quotient map $\chi\colon\Nsym\to\Sym$, sometimes referred to as the \emph{forgetful map} in that it allows image elements to commute.
In view of Equation \eqref{eqn:mono},   $\chi \left(M^*_\alpha\right) = h_{\partitionof{\alpha}}$.
Appropriately,  the $M^*_\alpha$ form a basis of $\Nsym$ called the \emph{noncommutative complete symmetric functions} \cite{GKLLRT}, and adopting the convention of \cite{bergeron-zabrocki-2} we denote these by $\bh_\alpha:=M^*_\alpha$.
Likewise, in view of Equation \eqref{eqn:schur}, for the basis in $\Nsym$ dual to the quasisymmetric Schur functions we have
\begin{equation} \label{eqn:qs-s-map}
\chi \left(\nqs_\alpha \right) \;=\; s_{\partitionof{\alpha}}.
\end{equation}

In this paper we shall refer to the basis $\{\nqs_\alpha\}$ of $\Nsym$ as the \emph{noncommutative Schur functions}, since they are dual to the quasisymmetric Schur functions and they are pre-images of the Schur functions under the forgetful map $\chi$.

\section{A noncommutative Littlewood-Richardson rule} \label{sec:NSymLR}

As in \cite{fulton-1}, combinatorial formulas for the skew Schur function indexed by $\nu/\lambda$  are obtained from Equation \eqref{eqn:schur} by simply extending the formula to tableaux of skew shapes:
\begin{equation} \label{eqn:skewSchurComb}
 s_{\nu/\mu} = \sum_{T\in SSRT(\nu/\mu)} \xx^{cont(T)} \quad
 = \sum_{T\in SRT(\nu/\mu)} L_{Des(T)}.
\end{equation}
It follows that $s_{\nu/\mu} \neq 0$ if and only if $\mu\subseteq\nu$.
We derive analogous formulas for skew quasisymmetric Schur functions.

\begin{proposition}[Combinatorial formulas for skew quasisymmetric Schur functions]
\label{prop:qschur-skew}
\begin{equation} \label{eqn:qschur-stdtabs-skew}
\qs_{\gamma\cskew\beta} = \sum_{T\in SSCT(\gamma\cskew\beta)} \xx^ {cont(T)} \quad = \sum_{T\in SCT(\gamma\cskew\beta)} L_{Des(T)}.
\end{equation}
\end{proposition}
\begin{proof}
We work from the SCT based formula of Equation \eqref{eqn:qschur} and
the coproduct rule for the fundamental basis of Equation \eqref{eqn:qscoprod}.
We first introduce some notation.
Given a SCT $T$ having $n$ cells, and an integer $k$, $0\leq k \leq n$,
we denote by $\mho_k(T)$ the SCT comprising the cells of $T$ with entries $\{1,\ldots,k\}$, the higher-numbered cells being added to the base shape. 
Similarly, we denote by $\Omega_k(T)$ the standardization of the SCT formed by removing from $T$ the cells numbered $\{1,\ldots,n-k\}$.
For example,

\[
T = \; \tableau{ 6& 4&1\\ 7&3 \\ \gst&\gst&5&2\\ \gst&8\\ \gst&\gst&9\\ \\}
\qquad  \mho_4(T) = \; \tableau{ \gst&4&1\\ \gst&3 \\ \gst&\gst& \gst&2\\ \gst& \gst\\ \gst&\gst& \gst}
\qquad  \Omega_5(T) = \; \tableau{ 2 \\ 3 \\ \gst&\gst&1\\ \gst&4\\ \gst&\gst&5\\}.
\]

We denote by $T+k$ the tableau obtained by adding $k$ to every entry of $T$.
If $S$ is a filling of the diagram  $\gamma \cskew\beta$
and $T$  a filling of the diagram $\beta \cskew \alpha $,
then we denote by $S\cup T$ the natural filling of the diagram of shape  $\gamma\cskew \alpha $.
Thus if $T$ is a standard  tableau having $n$ cells, and $k$ is an integer, $0\leq k \leq n$, we have \[ T = (\Omega_{n-k}(T)+k) \cup \mho_{k}(T). \]

When expanding $\Delta \qs_\gamma$ in terms of the fundamental basis using Equation \eqref{eqn:qschur}, for each $T\in SCT(\gamma)$ we obtain a summand $\Delta L_\delta$ where $\delta=Des(T)$.
Assuming $\gamma\vDash n$,  this summand in turn expands as
\[  \Delta L_\delta = \sum_{i=0}^n L_{\alpha_i}\otimes L_{\beta_i}, \]
where $|\alpha_i|=i$, $|\beta_i|=n-i$, and either $\alpha_i \beta_i = \delta$ or $\alpha_i \nearcat\beta_i = \delta$.
We observe that $\alpha_i=Des(\mho_i(T))$ and $\beta_i=Des(\Omega_{n-i}(T))$.

Let $T_l\in SCT(\beta)$  and
 $T_u\in SCT(\gamma\cskew \beta)$. Assuming $\gamma\vDash n$ and $\beta\vDash(n-m)$,
it is clear that $T=(T_l + m)\cup T_u\in SCT(\gamma)$,  with $\mho_m(T)=T_u$ and $\Omega_{n-m}(T)=T_l$.
Moreover, every SCT of shape $\gamma$ that splits in the above fashion into upper and lower parts with the lower part of shape $\beta$ is obtained in this way.
This establishes the formula
\[ \qs_{\gamma\cskew\beta}  = \sum_{T\in SCT(\gamma\cskew\beta)} L_{Des(T)}. \]
The remaining formula follows from Proposition  \ref{prop:standardization}:
\begin{eqnarray*}
\sum_{\widehat{T}\in SCT(\gamma\cskew\beta)} L_{Des(\widehat{T})}
&=& \sum_{\widehat{T}\in SCT(\gamma\cskew\beta)} \sum_{\tau\refines Des(\widehat{T})} \xx^ \tau \\ \\
&=& \sum_{\widehat{T}\in SCT(\gamma\cskew\beta)} \sum_{std(T)= \widehat{T}} \xx^{cont(T)}
\;=\; \sum_{T\in SSCT(\gamma\cskew\beta)}  \xx^{cont(T)}.
\end{eqnarray*}
\end{proof}

As a corollary of Proposition \ref{prop:qschur-skew}, it follows that  $\qs_{\gamma\cskew\beta} \neq0$ if and only if $\beta \leq_C\gamma$.

\begin{remark}
The skew quasisymmetric Schur functions generalize the classical skew Schur functions.
In particular, every  skew Schur function is equal to some skew quasisymmetric Schur function.
See Section  \ref{subsec:sym-SQS} for details.
\end{remark}

The structure constants of the Schur functions are called the Littlewood-Richardson coefficients.
The Littlewood-Richardson rule provides a combinatorial interpretation of these coefficients, one proof that they are nonnegative integers.
We state the rule here in terms of reverse tableaux.
Given a partition $\lambda$, we define the \emph{canonical reverse tableau of shape $\lambda$} to be  $\tilde{U}_\lambda := \rho(U_{\reverse{\lambda}})$.
Define $V_\lambda$ to be the unique SSRT of partition shape $\lambda$ and content $\reverse{\lambda}$, i.e., all 1's in the last row, 2's in the second to last row, etc.
Note that $std(V_\lambda)=\tilde{U}_\lambda$.
A  \emph{reverse Littlewood-Richardson tableau}  is an SSRT which rectifies to $V_\lambda$ for some partition $\lambda$.

\begin{example}
\[ \tilde{U}_{3221} \;\;=\;\; \tableau{8&7&6\\5&4\\3&2\\1 }
\qquad \qquad
 V_{3221} \;\;=\;\; \tableau{4&4&4\\3&3\\2&2\\1 }
 \]
\end{example}

\begin{theorem}[Classical Littlewood-Richardson rule] \label{thm:classical-LR}
In the expansion
\[ s_\lambda s_\mu = \sum_\nu c^\nu_{\lambda\mu} s_\nu, \quad\text{equivalently,}\quad
s_{\nu/\mu}  = \sum_\lambda c^\nu_{\lambda\mu} s_\lambda,  \]
the coefficient $c^\nu_{\lambda\mu}$ counts the number of reverse Littlewood-Richardson tableaux of shape $\nu/\mu$ and content $\reverse{\lambda}$, equivalently,
the number of SRT $T$ of shape $\nu/\mu $
such that $rect(T)=\tilde{U}_\lambda$.
\end{theorem}

A proof of Theorem \ref{thm:classical-LR} can be found in various texts, such as \cite{fulton-1}.
Our main result is a direct analogue of the Littlewood-Richardson rule, showing that the $\{\nqs_\alpha\}$ structure constants $C^\gamma_{\alpha\beta}$, which we will call the \emph{noncommutative Littlewood-Richardson coefficients}, have a combinatorial interpretation and thus are nonnegative integers.

\begin{theorem}[Noncommutative Littlewood-Richardson rule]  \label{thm:qschur-skew-pos}
Let the coefficients $C^\gamma_{\alpha\beta}$ be defined, equivalently, by any of the following expansions
\begin{eqnarray} \label{eqn:nc-coeffs}
\nqs_\alpha \cdot \nqs_\beta &=& \sum_\gamma C^\gamma_{\alpha\beta}\, \nqs_\gamma \\
\Delta \qs_{\gamma}  &=& \sum_{\alpha,\beta} C^\gamma_{\alpha\beta}\,\qs_\alpha \otimes \qs_\beta \\
\qs_{\gamma\cskew\beta}  &=& \sum_\alpha C^\gamma_{\alpha\beta}\, \qs_\alpha \\
C^\gamma_{\alpha\beta} &=& \langle \qs _{\gamma\cskew \beta }, \nqs_\alpha \rangle = \langle \qs _{\gamma }, \nqs_\alpha \cdot \nqs_\beta \rangle.
\end{eqnarray}
Then $C^\gamma_{\alpha\beta}$ counts the number of SCT $T$ of shape $\gamma\cskew\beta$ such that $rect(T)=U_\alpha$.
\end{theorem}


\begin{example}
\[
\begin{array}{ccccc}
  \tableau{  \gst &\gst &\gst &3&2 \\ \gst&\gst &1  }   
\quad & \quad \tableau{  \gst &\gst &\gst &2 \\ \gst&\gst &3&1  }   
\quad & \quad \tableau{  3&1\\ \gst &\gst &\gst &2 \\ \gst&\gst   }   
\quad & \quad \tableau{  3&1\\ \gst &\gst &\gst  \\ \gst&\gst &2  }   
\quad & \quad \tableau{  1\\ \gst &\gst &\gst &3&2 \\ \gst&\gst   }   
\\ \\  \tableau{ 3\\ \gst &\gst &\gst &2 \\ \gst&\gst&1   }   
\quad & \quad \tableau{ 1\\ \gst &\gst &\gst &2 \\ \gst&\gst &3  }   
\quad & \quad \tableau{  1\\3&2\\ \gst &\gst &\gst  \\ \gst&\gst   }   
\quad & \quad \tableau{ 1\\3\\ \gst &\gst &\gst  \\ \gst&\gst&2  }   
\quad & \quad \tableau{  1\\3\\ \gst &\gst &\gst &2 \\ \gst&\gst   }   
\end{array}
\]
\[ \nqs_{12}\cdot \nqs_{32} = \nqs_{53} +\nqs_{44} +\nqs_{242} +\nqs_{233}+\nqs_{152}+2\,\nqs_{143} +\nqs_{1232} +\nqs_{1133} +\nqs_{1142}  \]
\end{example}

The proof of this theorem will require several intermediate results, so we postpone the proof to Section \ref{sec:NSymLRproof}.
In the remainder of this section we consider some consequences of the main result.

The following family of decompositions of the classical Littlewood-Richardson coefficients follows immediately from Equation  \eqref{eqn:qs-s-map}.
\begin{corollary} Let $\alpha$ and $\beta$ be compositions with $\lambda=\partitionof{\alpha}$ and $\mu=\partitionof{\beta}$, and let $\nu$ be a partition.
Then $\chi(\nqs_\alpha\cdot \nqs_\beta) = s_\lambda\cdot s_\mu$, and
\begin{equation}
c^\nu_{\lambda\mu} = \sum_{\partitionof{\gamma}=\nu} C^\gamma_{\alpha\beta}.
\end{equation}
\end{corollary}

As special cases of Theorem  \ref{thm:qschur-skew-pos} we have Pieri rule analogues, that is, certain products in which the Littlewood-Richardson coefficients are all 0 or 1. 

\begin{corollary}[Noncommutative Pieri rules]
We have
\[ \nqs_{(n)}\cdot \nqs_\beta = \sum_\gamma \nqs_\gamma, \]  
where $\gamma$ runs over all compositions $\gamma \geq_C \beta$ such that $|\gamma\cskew\beta|=n$ and $\gamma\cskew\beta$ is a horizontal strip where the cells have been added from left to right.
Similarly,
\[ \nqs_{(1^n)}\cdot \nqs_\beta = \sum_\gamma \nqs_\gamma, \]  
where $\gamma$ runs over all compositions $\gamma \geq_C \beta$ such that $|\gamma\cskew\beta|=n$ and $\gamma\cskew\beta$ is a vertical strip where the cells have been added from right to left.
\end{corollary}

\hilight{One may compare this corollary with the classical Pieri rules for Schur functions \cite{fulton-1} 
and the Pieri rules for quasisymmetric Schur functions \cite[Theorem 6.3]{HLMvW-1}.}

%
%
\section{Proof of the noncommutative Littlewood-Richardson rule} \label{sec:NSymLRproof}

This section is devoted to proving Theorem \ref{thm:qschur-skew-pos}, by means of Proposition \ref{prop:C-complete} below.
In order to first outline the proof, we note the following view of Theorem \ref{thm:classical-LR} in relation to Equation \eqref{eqn:skewSchurComb}.
Haiman \cite{haiman-dualeq} defines the notion of \emph{dual equivalent} tableaux.
We use an equivalent definition, and for our purposes we find it convenient to restrict our attention to SRT, whose reading words may be viewed as permutations in one-line format.
Thus two SRT $T$ and $T'$ are dual equivalent if they have the same skew shape and $w_{col}(T) \stackrel{Q}{\sim} w_{col}(T')$. 
Note that this implies that $\pshape(T) =\pshape(T')$.
A restatement of \cite[Theorem 2.13]{haiman-dualeq} in our context is that 
each dual equivalence class of tableaux is ``complete'' 
in the sense that
if   $w$ is a permutation such that $w \stackrel{Q}{\sim} w_{col}(T)$
then there is a tableau $T'$ dual equivalent to $T$ such that $w_{col}(T')=w$. 
This implies on the one hand that for any dual equivalence class $[T]$ of SRT,
\[ \{\rho(rect(T')) : T'\in [T] \} = SRT(\lambda), \]
where $\lambda=\pshape(T)$, and thus \[ s_\lambda = \sum_{T'\in [T]} L_{Des(T')}. \]
On the other hand, this completeness implies that the set
$$\{T\in SRT \colon \rho(rect(T))=\tilde{U}_\lambda \text{ for some partition } \lambda \}$$ is a transversal (i.e., a set of representatives) of the collection of dual equivalence classes of SRT, or equivalently, that the set of  reverse Littlewood-Richardson tableaux forms a transversal of the collection of dual equivalence classes of SSRT.
These two familiar facts are embodied in Theorem \ref{thm:classical-LR}. 
The essence of our proof is to show that an analogous pattern holds for SCT.

\begin{definition}[$C$-equivalence] \label{def:C-equiv}
Permutations $\omega$ and $\pi$ are \emph{$C$-equivalent}, denoted $\omega\cequiv\pi$, if $\omega\stackrel{Q}{\sim}\pi$  and $\cshape(P(\omega))=\cshape(P(\pi))$.
We denote the $C$-equivalence class of the permutation $\pi$ by $[\pi]_C$.
The \emph{rectified shape} of $[\pi]_C$ is $\cshape(P(\pi))$.

Two SCT $T$ and $T'$  are $C$-equivalent $T\cequiv T'$ if they have the same skew shape and $w_{col}(T)\cequiv w_{col}(T')$.
We denote the $C$-equivalence class of $T$ by $[T]_C$.
The \emph{rectified shape} of $[T]_C$ is $\cshape(T)$.
We say that $[T]_C$ is \emph{complete} if
\begin{equation} \label{eq:C-eqs-T-complete}
\{ w_{col}(T') : T'\in [T]_C \} = [w_{col}(T)]_C.
\end{equation}
\end{definition}

\begin{example} Consider the permutations
$$3421\stackrel{Q}{\sim}2431\stackrel{Q}{\sim}1432$$for which
$$3421\not\cequiv 2431\cequiv 1432$$as $\rho^{-1}(P(\omega))$ for each $\omega$ is, respectively,
$$\tableau{3&2&1\\4}\qquad\tableau{2\\4&3&1}\qquad\tableau{1\\4&3&2}$$and hence $\cshape(P(\omega))$ for each $\omega$ is, respectively, $(3,1), (1,3), (1,3)$.
\end{example}

\begin{proposition} \label{prop:C-eqs}
Let $[\pi]_C$ be a $C$-equivalence class of permutations having rectified shape $\alpha$.
Then
\begin{equation} \label{eq:qschur-C-eqs-pi}
\qs_\alpha =  \sum_{\sigma\in[\pi]_C} L_{Des(P(\sigma))}.
\end{equation}
\end{proposition}
\begin{proof}
Let $f^\lambda=|SRT(\lambda)|.$
The RSK correspondence implies that there are $(f^\lambda)^2$ permutations $\pi$ such that $P(\pi)\in SRT(\lambda)$, and that these can be partitioned into $f^\lambda$ $Q$-equivalence classes, each class containing exactly one permutation from each $P$-equivalence class.
Therefore the rectifications of the respective permutations of a given $C$-equivalence class, say having rectified shape $\alpha$, form the complete set of SCT of shape $\alpha$.
\end{proof}

Note that the set of permutations
$\{ \pi : \rho^{-1}(P(\pi))=U_\alpha \text{ for some } \alpha \}$
is a transversal for the collection of $C$-equivalence classes of permutations.

\begin{proposition} \label{prop:C-complete}
$[T]_C$ is complete for every SCT $T$.
\end{proposition}

Proposition \ref{prop:C-complete} together with Equation  \eqref{eq:C-eqs-T-complete} imply  $\{ T \in SCT\colon rect(T)=U_\alpha \text{ for some } \alpha \}$
is a transversal for the collection of SCT $C$-equivalence classes,
and in conjunction with Propositions \ref{prop:qschur-skew} and  \ref{prop:C-eqs}
this is sufficient to prove Theorem \ref{thm:qschur-skew-pos}.
Before outlining the proof of  Proposition \ref{prop:C-complete},
we review some more known material \cite{assaf-1,haiman-dualeq,sagan} regarding Knuth and dual Knuth equivalence.

\medskip
Knuth and dual Knuth equivalence of permutations can be characterized by word transformations, or \emph{moves}.
If in the one-line notation of the permutation $\omega$ the elements $x=\omega_{k}$, $y=\omega_{k+1}$, and $z=\omega_{k+2}$ are not in monotonic order, then the \emph{elementary Knuth move} $p_k$ can be applied to $\omega$ by exchanging an appropriate pair of adjacent elements to obtain a Knuth equivalent permutation, specifically
\begin{equation}
 \omega = \ldots xyz \ldots \quad\Longrightarrow\quad
p_k(\omega)=\begin{cases}
\;\ldots yxz \ldots & \text{if } x<z<y\text{ or } y<z<x, \\
\;\ldots xzy \ldots & \text{if } y<x<z\text{ or } z<x<y .
\end{cases}
\end{equation}
Similarly, if 
the elements $\{k,k+1,k+2\}$, labeled $x$, $y$, and $z$ as they appear left to right in~$\omega$, are not in monotonic order (i.e., $y\neq k+1$), then the \emph{elementary dual Knuth move} $q_k$ can be applied to $\omega$ by exchanging $x$ and $z$ to obtain a dual Knuth equivalent permutation, i.e., \begin{equation}
  \omega = \ldots x \ldots y \ldots z \ldots \quad \Longrightarrow\quad
q_k(\omega)=\;\ldots z \ldots y \ldots x \ldots
\end{equation}

Knuth equivalence can be described as the transitive closure of the $p_k$ moves
while dual Knuth equivalence can be described as the transitive closure of the $q_k$ moves.
It is straightforward to verify that if $p_i$ and $q_j$ are both applicable to $\omega$, then the operators commute, that is,
\begin{equation} \label{eqn:pq-commute}
  p_i(q_j(\omega)) \;=\; q_j(p_i(\omega)).
\end{equation}
Among other consequences of these characterizations, if permutations are viewed as the column words of tableaux, then the $p_k$ moves preserve the descent sets of tableaux, and hence rectification preserves descents of tableaux, i.e., $Des(T) =Des(rect(T))$.

Given a partition $\lambda\vdash n$, let $H^\lambda$ be the undirected  graph whose vertex set $V(H^\lambda)$ is the set of all permutations $\pi\in S_n$ such that $sh(P(\pi))=\lambda$,
and where there is an edge $(\sigma,\pi)$  if and only if $\sigma$ and $\pi$ are related by an elementary dual Knuth move, i.e., $\sigma=q_k(\pi)$ for some $k$.
By definition,  the vertex sets of the connected components of  $H^\lambda$ are, respectively, the dual Knuth equivalence classes, and these classes are indexed by $SRT(\lambda)$.
Given $T\in SRT(\lambda)$, let $G^T$ be the connected component of $H^\lambda$ whose vertex set is indexed by $T$, i.e., $Q(\pi)=T$ for all $\pi\in V(G^T)$.
If $\sigma=p_k(\pi)$, and $Q(\sigma)=T$ and $Q(\pi)=T'$, then the relations of Equation \eqref{eqn:pq-commute} imply that $p_k$ defines a $P$-class preserving graph isomorphism between $G^T$ and $G^{T'}$.

\begin{proof}[Proof outline of Proposition \ref{prop:C-complete}]

The elementary moves  described above that define dual Knuth relations between permutations are of two types \cite{sagan}.
An elementary dual Knuth move of the first kind, denoted $\pi \stackrel{1*}{\cong} \sigma$ is of the form
\[ \pi = \ldots (k+1)\ldots k \ldots (k+2) \ldots \qquad \text{ and } \qquad \sigma = \ldots (k+2)\ldots k \ldots (k+1) \ldots, \]
whereas one of the second kind, denoted $\pi \stackrel{2*}{\cong} \sigma$ is of the form
\[ \pi = \ldots k \ldots (k+2) \ldots (k+1)\ldots \qquad \text{ and } \qquad \sigma = \ldots (k+1)\ldots (k+2)\ldots k  \ldots. \]
The descent set of a permutation (not to be confused with the descent set of a tableau) is defined as
\[ descents(\pi) \;:=\; \{ i : \pi_i > \pi_{i+1} \}. \]
Note that the dual Knuth moves preserve the descent set of the permutations, so $Q$-equivalent permutations have the same descent set.

Every $T\in SCT(\gamma\cskew\beta)$ determines a column word $\omega =w_{col}(T)$.
Consider the multiset of column numbers of the cells of $T$, written as a weakly increasing sequence $f$ of length $n=|\gamma\cskew\beta|$.
Writing $f$ as a word in this way, we see that
at those positions where there is a descent in $\omega$,
$f$ is strictly increasing.
We say that an arbitrary weakly increasing word is \emph{compatible} with $\omega$ if this condition holds.
Since $Q$-equivalent permutations have the same descent set, they share the same set of compatible words.

Clearly $colseq(T)$ can be recovered from the pair $(\omega,f)$, and hence $T$ itself can be recovered from the pair along with the base shape $\beta$, which we think of as
 ``applying'' $(\omega,f)$ to the base shape $\beta$.

\begin{example}
\[
\begin{array}{c}  \vspace{6pt}
\tableau{ 3 &1\\ \gst&\gst&\gst  \\ \gst&\gst &2\\ \gst&\gst&\gst&\gst &4\\ \gst&5 }
\\  T
\end{array}
\qquad \rightsquigarrow \qquad
\begin{array}{c}
\omega = w_{col}(T) = (31524)\\[5pt]
f = (12235)\\[5pt]
\beta=(3,2,4,1) \\[10pt]
T = (\omega , f)\beta
\end{array}
\]
\end{example}

Suppose that $T\in SCT(\gamma\cskew\beta)$ determines the pair $(\omega,f)$, i.e., $T= (\omega,f)\beta$, and that $\sigma\cequiv\omega$.
Then we want to show that  applying 
$(\sigma,f)$ to $\beta$ is defined, yielding $T'=(\sigma,f)\beta\in SCT(\gamma\cskew\beta)$. 
By Proposition \ref{prop:rho-general}, $T$ can be paired with a unique $\hat{T}\in SRT(\partitionof{\gamma}/\partitionof{\beta})$ such that $\omega=w_{col}(T)=w_{col}(\hat{T})$.
Now $\sigma\cequiv\omega$ implies $\sigma\stackrel{\tiny Q}{\sim}\omega$, and so by the completeness of dual equivalence classes of SRT, there exists $\hat{T}'\in SRT(\partitionof{\gamma}/\partitionof{\beta})$ such that $\sigma=w_{col}(\hat{T})$.
Again by Proposition \ref{prop:rho-general}, $\hat{T}'$ can be paired with a unique $T'\in SCT(\eta\cskew\beta)$ such that $\sigma=w_{col}(T')$ and $\partitionof{\eta}=\partitionof{\gamma}$,
and it follows that $T'=(\sigma,f)\beta$.
The only remaining question is whether $\eta=\gamma$.

The first step of our proof is to show that the proposition holds when $\sigma$ and $\omega$ are related by an elementary dual Knuth relation.
This step is further divided into two cases.
The
critical case comprises those relations of the second kind where the value $(k+2)$ appears in the first column of the skew tableau and the value $k$ (resp. $(k+1)$) appears in the second column of $T$ (resp. $T'$).
We will refer to this as the \emph{rigid case} and denote it by $T\Bumpeq T'$ (formally defined below).
We will start by considering first the easier situation which
covers all of the remaining cases.

Given a $C$-equivalence class of permutations $[\pi]_C$, with $Q(\pi)=U$ and rectified shape $\alpha$, let $G^U_\alpha$ be the subgraph of $G^U$  induced by the vertex set $V(G^U_\alpha)=[\pi]_C$.
Our final step is to show that $G^U_\alpha$ is always connected.
Connectivity of $G^U_\alpha$ implies that any two elements of $[\pi]_C$ can be transformed one into the other by a sequence of elementary dual Knuth moves.
This combined with the previous steps proves the proposition.
The steps of this outline are proved in the remainder of this section.
\end{proof}

\begin{definition}
Let $T\in SCT(\gamma\cskew\beta)$, $colseq(T)=(c_n,\ldots,c_1)$,
determining the pair $(\omega,f)$.  Let $\omega \stackrel{2*}{\cong} \sigma$, $T'=(\sigma,f)\beta$.
We say that $T$ and $T'$ are \emph{rigidly related}, denoted $T\Bumpeq T'$, if both
\begin{enumerate}
\item   $\omega$ is obtained from $\sigma$ by exchanging the values  $k$ and $k+1$, i.e., $\sigma=q_k(\omega)$.
\item  $c_{k+2}=1$, and $\{c_k,c_{k+1}\} = \{1,2\}$ as sets.
\end{enumerate}
\end{definition}

Note that the definition of $T\Bumpeq T'$ does \emph{not} assume that $T\cequiv T'$.
\begin{example}
\[ \begin{array}{cccl}
T \quad \tableau{ 3 & 1 \\ 5 & 4 \\ 6 \\ \gst & \gst &  \gst & 2 \\ \gst & 7 }
& \Bumpeq &
\tableau{ 3 & 1 \\ 4  \\ 6 &5 \\ \gst & \gst &  \gst & 2 \\ \gst & 7 }
 T' & k=4\\ \\
colseq(T)=(2112142) & & colseq(T')=(2121142) \\
(\omega,f)=(3561472,1112224) &&  (\sigma,f)=(3461572,1112224)
\end{array}
\]
\end{example}

\subsection{The easy case}

\begin{proposition}[First case] \label{prop:Qshape-case1}
Suppose $T\in SCT(\gamma\cskew\beta)$, determining the pair $(\omega,f)$.  Let $\sigma$ be a permutation such that either $\sigma \stackrel{1*}{\cong}\omega$ or $\sigma \stackrel{2*}{\cong}\omega$ 
and $T'=(\sigma,f)\beta \not\Bumpeq T$.
Then $T'\in SCT(\gamma\cskew\beta)$.
\end{proposition}
\begin{proof}
First, suppose that $\omega \stackrel{1*}{\cong} \sigma$, say
\[ \omega = \ldots (k+1)\ldots k \ldots (k+2) \ldots \qquad \text{ and } \qquad \sigma = \ldots (k+2)\ldots k \ldots (k+1) \ldots . \]
Suppose $(k+1)$ lies in column $i$ in $T$.
There must be a descent in $\omega$ somewhere between $(k+1)$ and $k$, so $k$ and $(k+2)$ lie strictly to the right of $(k+1)$ in $T$, say in columns $j$ and $j'$ respectively, when $i<j\leq j'$.
If $j'-i>1$, then all relationships between cells remain the same, that is, $T'$ is obtained from $T$ by simply exchanging places of the entries $(k+1)$ and $(k+2)$, and so $T$ and $T'$ have the same shape.
Otherwise we have $j'=j=i+1$, when $T'$ is obtained from $T$ by some permutation of the entries $k$, $(k+1)$, and $(k+2)$, so that $T$ and $T'$ again have the same shape:
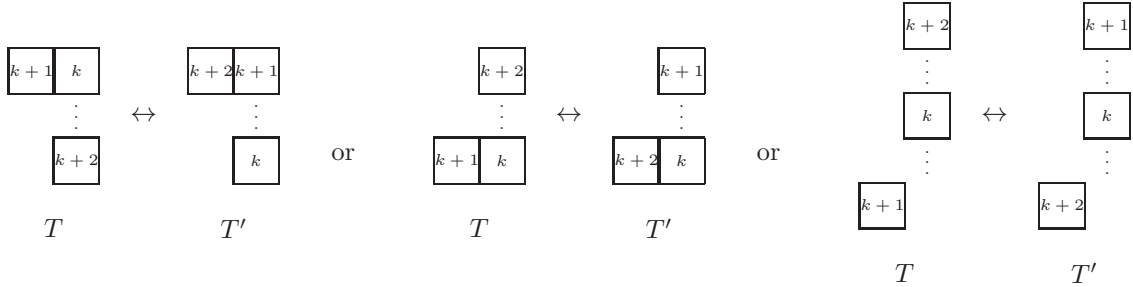
\begin{center}
\begin{tabular}{ccccc}

\setlength{\unitlength}{17pt}
\begin{picture}(6,7)(0,-3)  
\put(0,2){ \framebox(1,1){\tiny $k+1$} }
\put(1,2){ \framebox(1,1){\tiny $k$} }
\put(1,1.15){ \makebox(1,1){\tiny $\vdots$} }
\put(1,0){ \framebox(1,1){\tiny $k+2$} }
\put(0.5,-1.5){ \makebox(1,1){$T$} }
\put(2.5,1){ \makebox(1,1){$\leftrightarrow$} }
\put(4,2){ \framebox(1,1){\tiny $k+2$} }
\put(5,2){ \framebox(1,1){\tiny $k+1$} }
\put(5,1.15){ \makebox(1,1){\tiny $\vdots$} }
\put(5,0){ \framebox(1,1){\tiny $k$} }
\put(4.5,-1.5){ \makebox(1,1){$T'$} }
\end{picture}

& \raisebox{60pt}{ $\quad$or$\quad$ } &

\setlength{\unitlength}{17pt}
\begin{picture}(6,7)(0,-3)  
\put(1,2){ \framebox(1,1){\tiny $k+2$} }
\put(1,1.15){ \makebox(1,1){\tiny $\vdots$} }
\put(0,0){ \framebox(1,1){\tiny $k+1$} }
\put(1,0){ \framebox(1,1){\tiny $k$} }
\put(0.5,-1.5){ \makebox(1,1){$T$} }
\put(2.5,1){ \makebox(1,1){$\leftrightarrow$} }
\put(5,2){ \framebox(1,1){\tiny $k+1$} }
\put(5,1.15){ \makebox(1,1){\tiny $\vdots$} }
\put(4,0){ \framebox(1,1){\tiny $k+2$} }
\put(5,0){ \framebox(1,1){\tiny $k$} }
\put(4.5,-1.5){ \makebox(1,1){$T'$} }
\end{picture}

& \raisebox{60pt}{ $\quad$or$\quad$ } &

\setlength{\unitlength}{17pt}
\begin{picture}(6,7)(0,-2)  
\put(1,4){ \framebox(1,1){\tiny $k+2$} }
\put(1,3.15){ \makebox(1,1){\tiny $\vdots$} }
\put(1,2){ \framebox(1,1){\tiny $k$} }
\put(1,1.15){ \makebox(1,1){\tiny $\vdots$} }
\put(0,0){ \framebox(1,1){\tiny $k+1$} }
\put(0.5,-1.5){ \makebox(1,1){$T$} }
\put(2.5,2){ \makebox(1,1){$\leftrightarrow$} }
\put(5,4){ \framebox(1,1){\tiny $k+1$} }
\put(5,3.15){ \makebox(1,1){\tiny $\vdots$} }
\put(5,2){ \framebox(1,1){\tiny $k$} }
\put(5,1.15){ \makebox(1,1){\tiny $\vdots$} }
\put(4,0){ \framebox(1,1){\tiny $k+2$} }
\put(4.5,-1.5){ \makebox(1,1){$T'$} }
\end{picture}

\end{tabular}
\end{center}
Thus in the case $\omega \stackrel{1*}{\cong} \sigma$, $T$ and $T'$ have the same shape.

Next consider the case $\omega \stackrel{2*}{\cong} \sigma$, say
\[ \omega = \ldots k \ldots (k+2) \ldots (k+1)\ldots \qquad \text{ and } \qquad \sigma = \ldots (k+1)\ldots (k+2)\ldots k  \ldots . \]
Suppose $(k+1)$ lies in column $j$ in $T$.
There must be a descent in $\omega$ somewhere between $(k+2)$ and $(k+1)$, so $k$ and $(k+2)$ lie strictly to the left of $(k+1)$ in $T$, say in columns $i$ and $i'$ respectively, when $i\leq i'< j$.
Again, if $j-i>1$, then all relationships between cells remain the same, that is, $T'$ is obtained from $T$ by simply exchanging places of the entries $k$ and $(k+1)$, and so $T$ and $T'$ have the same shape.
Otherwise we have $i=i'=j-1$.
Thus, as we construct $T$ by applying $(\omega,f)$ to $\beta$, cells $(k+2)$ and $k$ are added to rows of length $i-1$.
Since we are assuming that $T\not\Bumpeq T'$,  we have $i>1$, and so $k$ will be inserted in a row below $(k+2)$.
Again, $T'$ is obtained from $T$ by simply exchanging places of the entries $k$ and $(k+1)$:
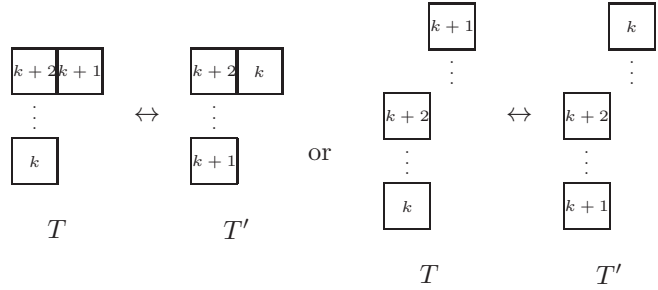
\begin{center}
\begin{tabular}{ccc}

\setlength{\unitlength}{17pt}
\begin{picture}(6,7)(0,-3)  
\put(0,2){ \framebox(1,1){\tiny $k+2$} }
\put(1,2){ \framebox(1,1){\tiny $k+1$} }
\put(0,1.15){ \makebox(1,1){\tiny $\vdots$} }
\put(0,0){ \framebox(1,1){\tiny $k$} }
\put(0.5,-1.5){ \makebox(1,1){$T$} }
\put(2.5,1){ \makebox(1,1){$\leftrightarrow$} }
\put(4,2){ \framebox(1,1){\tiny $k+2$} }
\put(5,2){ \framebox(1,1){\tiny $k$} }
\put(4,1.15){ \makebox(1,1){\tiny $\vdots$} }
\put(4,0){ \framebox(1,1){\tiny $k+1$} }
\put(4.5,-1.5){ \makebox(1,1){$T'$} }
\end{picture}

& \raisebox{60pt}{ ~~~~~or~~~~~ } &

\setlength{\unitlength}{17pt}
\begin{picture}(6,7)(0,-2)  
\put(1,4){ \framebox(1,1){\tiny $k+1$} }
\put(1,3.15){ \makebox(1,1){\tiny $\vdots$} }
\put(0,2){ \framebox(1,1){\tiny $k+2$} }
\put(0,1.15){ \makebox(1,1){\tiny $\vdots$} }
\put(0,0){ \framebox(1,1){\tiny $k$} }
\put(0.5,-1.5){ \makebox(1,1){$T$} }
\put(2.5,2){ \makebox(1,1){$\leftrightarrow$} }
\put(5,4){ \framebox(1,1){\tiny $k$} }
\put(5,3.15){ \makebox(1,1){\tiny $\vdots$} }
\put(4,2){ \framebox(1,1){\tiny $k+2$} }
\put(4,1.15){ \makebox(1,1){\tiny $\vdots$} }
\put(4,0){ \framebox(1,1){\tiny $k+1$} }
\put(4.5,-1.5){ \makebox(1,1){$T'$} }
\end{picture}

\end{tabular}
\end{center}
Thus in the case that $\omega \stackrel{2*}{\cong} \sigma$ and $T\not\Bumpeq T'$, $T$ and $T'$ have the same shape.
\end{proof}

\subsection{The rigid case}

\begin{proposition}[Second case] \label{prop:Qshape-case2}
Suppose $T\in SCT(\gamma\cskew\beta)$,  determining the pair $(\omega,f)$.  Let $\sigma$ be a permutation such that $\sigma\cequiv\omega$  and $T'=(\sigma,f)\beta \Bumpeq T$.
Then $T'\in SCT(\gamma\cskew\beta)$.
\end{proposition}

We prove Proposition \ref{prop:Qshape-case2} via its contrapositive.
Specifically, we assume that $T\Bumpeq T'$  \emph{without} assuming that $\cshape(P(\omega))=\cshape(P(\sigma))$.
Suppose $T\in SCT(\gamma\cskew\beta)$ and $T'\in SCT(\gamma'\cskew\beta)$.
We assume that $\gamma\neq\gamma'$, whence it suffices to show that  $\cshape(P(\omega))\neq\cshape(P(\sigma))$.

Without loss of generality, we assume that in $T$ the cells with $(k+1)$ and $(k+2)$ are in column 1,   in rows $r$ and $r+1$ respectively, that the cell with $k$ is in column 2,
and that  $\sigma=q_k(\omega)$, that is, $\sigma$ is obtained from $\omega$ by exchanging the values  $k$ and $(k+1)$.
The relations \eqref{eqn:pq-commute} then  imply that  $w_{col}(rect(T'))=q_k(w_{col}(rect(T)))$,
and hence $rect(T)\Bumpeq rect(T')$.  
We note that $T$ and $T'$ differ only in rows $r$ and $r+1$, the elements within these rows being re-arranged and all other rows being identical between $T$ and $T'$.
We will say that this is the pair of adjacent rows associated with the move $q_k$ on $T$,
or simply
 \emph{the row pair} of $T$ (or $T'$, when clear from context).
It follows from $\gamma\neq\gamma'$ and the configuration of entries $(k+2)$, $(k+1)$, and $k$ in $T$
that 
$\gamma_r >\gamma_{r+1}$, the respective row lengths necessarily being reversed in $T'$.
Our idea is to show that
the rows in
the row pair associated with $q_k$ on $rect(T)$ are also of different lengths, the lengths being reversed when comparing  $rect(T)$ to $rect(T')$.

For $1\leq j\leq\gamma_{r+1}$ we will refer to the configuration of three cells $(r,j)$, $(r,j+1)$, and $(r+1,j)$ as \emph{the $j$-th triple} of the row pair of $T$.
We say that the $j$-th triple is \emph{rigid} if either
\begin{enumerate}
\item  $j=1$ and the set of the entries in the triple is $\{k,k+1,k+2\}$ for some $k$, or
\item if $T(r+1,j)<T(r,j+1)$.
\end{enumerate}
We say that \emph{the row pair is rigid} if the two row lengths differ and all of its triples are rigid.
These conditions imply that $T(r+1,j)<T(r,j)$ for $j>1$ and that row $r$ is longer than row $r+1$.

In our context, the row pair of $T$ must be rigid, for if not, say if $j$ is the least index for which the triple is not rigid, then the entries in the cells in the $(j+1)$-st column of the row pair of $T$ and $T'$, as well as all columns to the right, would be the same between $T$ and $T'$, and hence $\gamma=\gamma'$, contrary to assumption.
In the two examples below, the row pair consists of the second and third row, with $k=7$.  In the first example the second triple is not rigid, while in the second example the row pair is rigid.
\[
\begin{array}{ccc|ccc}
\tableau{ 1  \\ 8 & 7 & \gxS{4} &3\\ 9 &5 &\gxS{2} \\ \gstD & \gstD & 10 & 6 }
&\quad \stackrel{k=7}{\longleftrightarrow} \quad &
\tableau{ 1  \\ 7 & 5 & \gxS{4} &3 \\ 9 &8 &\gxS{2} \\ \gstD & \gstD & 10 & 6 }
\qquad & \qquad
\tableau{ 3  \\ 8 & 7 & 5 &2 \\ 9 &4 &1 \\ \gstD & \gstD & 10 & 6 }
&\quad \stackrel{k=7}{\longleftrightarrow} \quad &
\tableau{ 3  \\ 7 & 4 &1 \\ 9 &8 &5&2 \\ \gstD & \gstD & 10 & 6 }
\\  & & & & & \\
T &not\; rigid & T' & T & rigid & T'
\end{array}
\]

Conversely, let $r'$ and $r'+1$ be the corresponding row pair of $rect(T)$ containing the entries $k+1$ and $k+2$ respectively.
Since, as noted above, $rect(T)\Bumpeq rect(T')$ and $w_{col}(rect(T'))=q_k(w_{col}(rect(T)))$,
if this row pair in $rect(T)$ is rigid, then $rect(T)$ and $rect(T')$ also have  different shapes.
Before proceeding with the proof of Proposition \ref{prop:Qshape-case2}, we shall need some intermediate results.

\begin{comment} 
\[
\begin{array}{c} \\
{\tiny
\tableau{ 
\bas{1} & \bas{2} \\
k\!+\!1 &k \\
k\!+\!2 } }  \\ \\
\text{first triple }    \\
\text{first column: rigid}
\end{array}
\qquad
\begin{array}{c} \\
{\tiny \tableau{  \bas{j} & \bas{j\!+\!1} \\ c & a \\  b} } \\ \\
\text{$j$-th triple, } j >1 \\
  b < a  \text{ : rigid}
\end{array}
\qquad
\begin{array}{c} \\
{\tiny \tableau{  \bas{j} & \bas{j\!+\!1} \\ c & a \\  b} } \\ \\
\text{$j$-th triple, } j >1 \\
  a < b  \text{ : not rigid}
\end{array}
\]
\end{comment} 

\medskip
An alternative method for computing the straight SSRT  $\rho(rect(S))$ of a skew SSRT $S$ is provided by using Schensted insertion, 
i.e., by successively inserting the entries of $w_{col}(S)=w_1w_2\cdots w_n$ into an initially empty tableau, that is, \[ \rho(rect(S)) = ((\emptyset \gets w_1)\gets w_2)\cdots \gets w_n.  \]
We define insertion for SSCT $S$ via Mason's bijection,
viz.
\[ (S\gets k) \;:=\; \rho^{-1}\left(\, \rho(S)\gets k\,\right). \]
Thus an alternative method for computing the straight SSCT $rect(S)$ is to insert $w_{col}(S)$ into an initially empty composition tableau.
We recall from \cite{mason-1} an explicit description of the insertion algorithm for SSCT.
In this context we regard a composition diagram $\alpha$ as a subset of the rectangular $\ell\times (m+1)$ array of cells where $\ell=\ell(\alpha)$ and $m$ is the largest part of $\alpha$.
The \emph{scanning order} of cells in this array is down each successive column starting at the rightmost column.  That is, cell $(i,j)$ is scanned before cell $(i',j')$ if $j>j'$ or if $j=j'$ and $i<i'$.
To insert a new element $k$  into an SSCT $T$, we apply the following algorithm to the cells in scanning order.

\begin{algorithm}[SSCT insertion] \label{alg:insertion} 
To compute $(T\gets k)$,
\begin{enumerate}
\item  Initialize the variable
$\;z := k$.
\item  If we are in the first column, place $z$ in the first cell of a new row such that the entries in the first column are increasing top to bottom, and halt.
\item  If the current cell $(i,j)$ is empty, and $(i,j-1)$ is not empty, and $z\leq T(i,j-1)$, then place $z$ in position $(i,j)$ and halt.
\item  If $(i,j)$ is not empty, and $T(i,j)<z\leq T(i,j-1)$, then swap $z$ with the entry in $T(i,j)$ (we say that the entry in $T(i,j)$ is ``bumped'') and continue.
\item Go to step (2), processing the next cell in scanning order.
\end{enumerate}
\end{algorithm}

\begin{example} In these examples, the cells of the insertion path are highlighted.
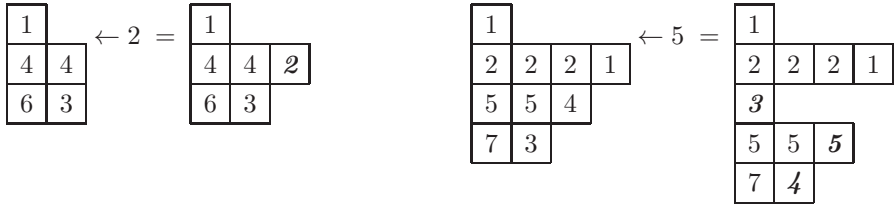
\[
\tableau{1\\4&4 \\ 6& 3 } \gets 2 \;=\; \tableau{1\\4&4&\gx{2} \\ 6& 3}
\qquad \qquad \qquad
\tableau{1\\2&2&2&1 \\ 5& 5&4 \\ 7&3} \gets 5 \;=\; \tableau{1\\2&2&2&1 \\ \gx{3}\\ 5& 5&\gx{5} \\ 7&\gx{4}}
\]
\end{example}
We note the following facts regarding insertion.
See \cite{fulton-1} and  \cite{HLMvW-1} for details.
\begin{itemize}
\item If a cell is on the insertion path, its contents are replaced with a larger value.
\item The entries of the cells of the insertion path in the new tableau in scanning order are strictly decreasing.
\item At most one cell per row is on the insertion path.
\item   If one successively inserts a strictly increasing sequence of elements into a reverse tableau, i.e., $T' = ((T\gets x_1) \gets x_2)\gets \cdots \gets x_k$ where $x_1<\cdots<x_k$, then the skew shape $\nu/\mu$, where $sh(T)=\mu$ and $sh(T')= \nu $, is a vertical strip.  This implies that under the bijection $\rho$, if $\cshape(T)=\beta$ and $\cshape(T')=\gamma$, then $\gamma\cskew\beta$ is also a vertical strip. 
\end{itemize}

\begin{proposition} \label{prop:rigid-rows}
Suppose SSCT $T$ has distinct entries and has a rigid row pair $r,r+1$.
Let
$\;U= (T\gets z)$ be the result of inserting an element $z$ into $T$, where we assume $z$ is not already an entry in $T$.
Then the corresponding row pair of $U$ is also rigid.
\end{proposition}
\begin{proof}
At most
one cell per row can be in the insertion path, and neither of them can lie in the first triple of the row pair.
If the insertion adds a new cell to the end of row $r$, then clearly the row pair in $U$ is also rigid.
If the insertion path does not contain any cell of the row pair, or if it contains a cell of row $r$ but not any cell of row $r+1$, then since any affected entry is replaced by a larger one, all the triples of the row pair in $U$ are also rigid, and so the row pair in $U$ is rigid.

Suppose that the insertion path contains the cells $(r,i)$ and $(r+1,j)$.
Now $i<j$ would imply that \[ U(r,i) < U(r+1,j)<U(r+1,i)=T(r+1,i) <T(r,i) <U(r,i), \]
a contradiction.
Also, $i=j$ would imply that \[ T(r+1,j-1)<T(r,j)=U(r+1,j)<U(r+1,j-1)=T(r+1,j-1), \]
again a contradiction.
Thus $j<i$ (where possibly $(r+1,j)$ is empty in $T$) and we have $U(r+1,j)<U(r,i)\leq U(r,j+1)$,
and so all triples of the row pair in $U$ are also rigid.

The remaining cases are when the insertion path contains the cell $(r+1,j)$, $j>1$, but no cell of row $r$,
in which case  $U(r+1,j) < T(r+1,j-1)<T(r,j)$.
Possibly $(r+1,j)$ is empty in $T$.
In any case, we need to show that $(r,j+1)$ is not empty in $T$ and that $U(r+1,j) <U(r,j+1)$.
Consider the value of the variable $z$ of the insertion algorithm at the point that it was processing the cell position $(r+1,j+1)$.
Since the cell $(r,j+1)$ is not on the insertion path, either $z<T(r,j+1)$ or $T(r,j)<z$.
In the former case, $U(r+1,j) <z<T(r,j+1)=U(r,j+1)$ and we are done.

In the latter case, take $T(r,j+1)$ to be 0 if $(r,j+1)$ is empty in $T$.
Suppose that $U(r+1,j) >T(r,j+1)$.
Then there must exist a cell $(s,i)$ on the insertion path lying strictly between $(r+1,j+1)$ and $(r+1,j)$ in scanning order such that
\[ T(r,j+1)<U(r+1,j)\leq T(s,i)   <T(r,j)< U(s,i)\leq z. \]
If $i=j+1$, which requires $r+1<s$, then $T(r,j)$, $T(r,j+1)$, and $T(s,j+1)$ would violate the definition of an SCT.
Otherwise $i=j$, which requires $s<r$, in which case we have
\[ T(s,j) <{T}(r,j) < U(s,j) <  U(s,j-1) =T(s,j-1) \]
and so $T(s,j-1)$, $T(s,j)$, and $T(r,j)$ would violate the definition of an SCT.
Thus our supposition is false; we must have $U(r+1,j) <T(r,j+1)$ as desired.
\end{proof}


To state the next proposition, we extend our notation for indexing
cells.
Let $X$ be a subset of the cell entries in the first column of a
tableau $T$.
We define $T(X)$ to be the set of those rows containing an element of
$X$ in its first column,
and we define $T(X,j)$ to be the set of  cell entries in the $j$-th  
column of the rows $T(X)$. For $X=\{r\}$ we also write
$T(\{r\},j)=x$,
omitting the
brackets on the right hand side.

\begin{example}[of notation]
\[
T = \quad \tableau{ 3 &1\\ 5 \\ 7 & 6 & 4 \\ \gst&\gst &2 }
\qquad
\begin{array}{rcl}
\\
T(\{{3},{7}\},2)&=&\{{1},{6}\} \\  T(\{{3},{7}\},3)&=&\{{4}\}\\
T(\{7\},3)&=&4
\end{array}
\]
\end{example}

\begin{proposition} \label{prop:T-and-P}
Let $T=(\omega,f)\beta$ be an SCT of shape $\gamma\cskew\beta$
with $\omega=w_{col}(T)=C_1\cdots C_t$ where $C_j$ is the set of  
entries in column $j$, and such that $C_1\neq\emptyset$.
Let $P_j$ be the partial rectification of $\omega$ obtained after  
inserting the prefix $C_1\cdots C_j$ of $\omega$  into the empty  
tableau.
Let $I_1 = C_1$ (as a set of cell entries) and use ${i}\in I_1$ to  
index rows of both $T$ and the partial rectifications.
Let $m = \max_{{i}\in I_1} \ell(row_{{i} }(T))$,
the maximum length over  $T(I_1)$.
Then for all $1\leq j \leq m$ we have the following. 
\begin{enumerate}
\item  The maximum row length in $P_j$ is $j$.
\item
Letting $I_j$ be the set of first column cell entries of those rows
  of $P_j$ of length $j$,
  $I_j$ also indexes the set of all rows in $T$
  that begin in column 1 and have length at least $j$.
\item  The entries $T(I_j,j) = P_j(I_j,j)$, and are in the same  
relative order within the column.
\end{enumerate}
\end{proposition}

\begin{example}[for Proposition  \ref{prop:T-and-P}]
\label{ex:T-and-P}
\[ T = \quad \tableau{ 3 &1\\ 11 &10 &8 &7 \\ 12 &6 &4 \\ \gst &13 &2
\\ \gst&\gst&\gst  &9\\ \gst&\gst&\gst  &5  } \]

\[ \begin{array}{c|c|c|c} \vspace{6pt}
\tableau{ \\ \gx{3} \\ \gx{11} \\ \gx{12} }
\quad & \quad
\tableau{ \\ 3 &\gx{1} \\ 11 &\gx{10} \\ 12 &\gx{6} \\ 13 }
\quad & \quad
\tableau{ \\ 3 &2 \\ 11 &10 &\gx{8} \\ 12 &6 &\gx{4} \\ 13 &1 }
\quad & \quad
\tableau{ 1 \\ 3 &2 \\ 4 \\ 11 &10 &9 & \gx{7} \\ 12 &8 &5 \\ 13 &6 }
\\
P_1 & P_2 & P_3 & P_4 \\
C_1= & C_1C_2= & C_1C_2C_3= & C_1C_2C_3C_4= \\
3,11,12 & 3,11,12,1,6,10,13 & 3,11,12,1,6,10,13, & 3,11,12,1,6,10,13,
\\
& &2,4,8 &
2,4,8,5,7,9 \\
I_1=\{{3},{11},{12}\} & I_2=\{{3},{11},{12}\} & I_3=\{{11},{12}\} &  
I_4=\{{11}\}
\end{array}\]

\medskip
In this example, $m=4$ and $rect(T)=P_4$.  The highlighted cells of  
$P_j$ match those of their counterparts in $T$.
\end{example}

\begin{proof}[Proof of Proposition \ref{prop:T-and-P}]
Proceed by induction on $j$.  The proposition clearly holds for $j=1$.
Hence we now assume that $j>1$ and that $C_j$ is nonempty.
When we insert $C_j$ into $P_{j-1}$ to obtain $P_j$,
we are adding a vertical strip to the overall shape of $P_{j-1}$ to  
obtain the shape of~$P_j$.
By hypothesis the longest rows of $P_{j-1}$ have length $j-1$, so it  
follows that $I_j\subseteq I_{j-1}$, and that the maximal row length
in $P_j$ is $j$, establishing part~(1). Now we prove the remaining
parts
together.

Let $I_j=\{r_1,\ldots, r_s\}$ with $r_1< \cdots < r_s$.
As $I_j\subseteq I_{j-1}$, by induction $T(\{r_i\},j-1)$ is nonempty,
say
$T(\{r_i\},j-1)=x_i$, for $i=1,\ldots,s$.
In fact, $T(\{r\},j-1)=P_{j-1}(\{r\},j-1)$ for all $r\in I_{j-1}$.
Then, by the definition of the insertion process of $C_j$ into
$P_{j-1}$,
$P_j( \{r_1\} ,j)=y_1$ where $y_1$ is the largest element of $C_j$  
such that $y_1<x_1$.
Observe that in particular, $\min C_j < x_1$ but
$\min C_j >P_{j-1}(\{r\},j-1)=T(\{r\},j-1)$ for all $r\in I_{j-1}$  
with $r<r_1$.
Note that during the insertion, it is possible that the entry $x_1$
in the
cell $(\{r_1\},j-1)$ might be replaced by an entry $x_1' > x_1$,
but since the elements of $C_j$ are inserted in increasing order,
this
must occur after $y_1$ has been inserted into the cell $(\{r_1\},j)$,
and these subsequent insertions do not affect the contents of  
$P_j(\{r_1\},j)$.

We claim that $T(\{r_1\},j)=y_1$.
As $y_1$ is the largest element in~$C_j$
that is smaller than $x_1$, by the triple condition for SCT
$y_1$ cannot be lower than $x_1$ in $T$.
If $y_1=T(\{r\},j)$ for some $r< r_1$,
then $r\in I_{j-1}$  and $T(\{r\},j-1)>y_1$,
a contradiction to the observation above.
Thus we have $y_1=T(\{r_1\},j)$.

Likewise, $P_j(\{r_2\},j)=y_2$ where  $y_2$ is the largest element of
$C_j\setnot\{y_1\}$ such that $y_2<x_2$, and similar reasoning as
above
yields $y_2=T(\{r_2\},j)$; continuing along these lines we find
$P_j(\{r_i\},j)=T(\{r_i\},j)$ for all $i\leq s$.

It remains to show that $T(\{r\},j)$ is empty for all
$r\in I_{j-1}\setminus I_j$.
Assume that $T(I_{j-1}\setminus I_j,j)$ is nonempty,
say $ \{r_{s+1},\ldots,r_t\} \subseteq I_{j-1}\setminus I_j$ with
$t>s$
are the additional row indices
with nonempty $T(\{r_k\},j)=y_k$;
we emphasize that the corresponding rows are not necessarily below
the row indexed $\{r_s\}$ but that these row indices may be
interleaved with
the ones in~$I_j$.
By induction, $T(\{r_k\},j-1)=x_k=P_{j-1}(\{r_k\},j-1)$ for all
$k\leq t$.
By definition of an SCT, $x_k>y_k$ for all $k\leq t$;
in particular, $P_{j-1}$ has at least $t$ rows that will be extended 
when the $t$ smallest elements of $C_j$ are inserted. 
Note that by definition of the bumping process, a box
which is filled at some step will never be vacated later,
the entry might only be replaced.
Now by definition of $I_j$, only the
$P_j(\{r_k\},j)$ for $k=1,\ldots,s$ are nonempty.
Hence we must have $s=t$, reaching the final contradiction.
This ends the proof of the case $j$
and the proposition now follows by induction.
\end{proof}


Now we are ready to complete the proof of Proposition \ref{prop:Qshape-case2}.
Set $X=\{k+1,k+2\}$.
Recall that the row pair $T(X)$  is rigid, and that we need to show that $rect(T)(X)$  is also rigid.
Let $m$ be one more than the length of the second (shorter) row of the pair  $T(X)$.
Using the notation of Proposition \ref{prop:T-and-P}, we claim that  $P_m(X)$ is rigid,
which by Proposition \ref{prop:rigid-rows} implies that  $rect(T)(X)$  is also rigid.
Proposition \ref{prop:T-and-P} implies that the row lengths of $P_m(X)$ are $m$ and $m-1$ respectively.
We show that the triples of $P_m(X)$ are rigid by induction on their column number.
The first (column 1) triple of $P_m (X)$ is the same as that of $T(X)$, which is rigid.
If $m=2$, then we are done.

\[
\begin{array}{c}
{\tiny \tableau{  k\!+\!1 &k & y & \bas{\cdots} \\  k\!+\!2 & x & z & \bas{\cdots} } } \\ \\ T
\end{array} \qquad \Longrightarrow \qquad
\begin{array}{c}
{\tiny \tableau{  k\!+\!1 &k \\  k\!+\!2 & x } } \\ \\ P_2
\end{array} \gets \quad C_3
\qquad \Longrightarrow \qquad
\begin{array}{c}
{\tiny \tableau{  k\!+\!1 &k & y \\  k\!+\!2 & x' & z } } \\ \\ P_3
\end{array}
\]
Otherwise, $m>2$.
Suppose that in $T$, $y$ is to the immediate right of $k$, and to the immediate right of $(k+2)$ we have elements $x$ and $z$, as shown in the diagram, where possibly $z$ could be empty.
By Proposition \ref{prop:T-and-P},  $x$ is to the immediate right of $(k+2)$ in $P_2$, and the third column of $P_3$ has $y$ and $z$ in the respective rows.
The element $x'$ to the immediate right of $(k+2)$ in $P_3$ could be $x$ or it could be some larger element due to bumping.
We claim that $x' < y$.
Since the second triple (column 2) of $T(X)$ is rigid, we have $x<y$.
During the insertion $P_3=P_2\gets C_3$,
once $y$ was in cell position $P_3(\{k+1\},3)$, any larger values inserted must have been larger that $k$ (else $y$ would have been bumped), and hence larger than $k+2$, leaving position $P_3(\{k+2\},2)$ unchanged.
Thus if $x$ was bumped during the insertion by $x'>x$, it had to have been bumped before or during placement of $y$ in position $P_3(\{k+1\},3)$, and this implies that $x' < y$.
Thus both the first and second triples of $P_3(X)$ are rigid.
Continuing the argument by induction, we obtain that all the triples of columns  $i$, $1\leq i < j$, of $P_j(X)$  are rigid for $1<j\leq m$,
and hence $P_m(X)$ is rigid as claimed.
This completes the proof of Proposition \ref{prop:Qshape-case2}.

\subsection{Connectivity of $G^U_\alpha$}

\begin{proposition} \label{prop:CT-graph}
Every graph  $G^U_\alpha$ is connected.
\end{proposition}
\begin{proof}
As noted above, if $U,U'\in SRT(\lambda)$ and permutations $\sigma\in V(G^U)$ and $\pi\in V(G^{U'})$ such that $\sigma=p_k(\pi)$ for some $k$, then $p_k$ defines a $P$-class preserving graph isomorphism between $G^U$ and $G^{U'}$, which therefore restricts to an isomorphism between $G^U_ \alpha $ and $G^{U'}_ \alpha$ for all $\partitionof{\alpha}=\lambda$.
Moreover,  just as $G^U$ is connected, the corresponding graph defined on $P$-equivalence classes using elementary Knuth moves for edges is also connected.
By transitivity, $G^U_\alpha$ is connected for all $U\in SRT(\lambda)$ if $G^U_\alpha$ is connected for any single $U\in SRT(\lambda)$.
For convenience we choose to work with the unique $U\in SRT(\lambda)$ such that vertices of $G^U$ are the respective column reading words of all $SRT(\lambda)$, namely that in which the elements are numbered consecutively in decreasing fashion down each column starting with the first column, for example, \[ \tableau{ 11 &7 &4 &1\\ 10 &6&3 \\ 9 &5&2 \\ 8 }.  \]
Accordingly we may drop $U$ from our notation and refer to $G^U_\alpha$ simply as $G_\alpha$.

Let $T$ be a straight SCT having $m$ columns, and let $C_j$ be the set of entries in the $j$-th column.
We now identify $T$ with its column word and its column tabloid via
\[   T \quad \leftrightarrow \quad w_{col}(T) \quad \leftrightarrow \quad(C_1,\ldots,C_m). \]
This means we can perform our analysis in terms of SCT  by identifying  $V(G_\alpha)$ with $SCT(\alpha)$.

An elementary  dual Knuth move applied to a permutation in one-line notation exchanges the positions of the elements $k$ and $k+1$ for some value $k$, subject to the condition that  either $k+2$ or $k-1$ is positioned somewhere between $k$ and $k+1$.
We say that the  move is \emph{applicable} to the permutation if it meets the required condition.
Applying a move to $w_{col}(T)$ is equivalent to exchanging two elements between different column sets in the tabloid.
Expressed in terms of the tableau itself,
the applicability condition is that we may exchange $k$ and $k+1$ as long as
\begin{enumerate}
\item  $k$ lies  strictly to the left of $k+1$ and
   \begin{enumerate}
   \item  $k+2$ lies  strictly to the left of $k+1$ and weakly right of $k$, or
   \item  $k-1$ lies  weakly to the left of $k+1$ and strictly right of $k$.
   \end{enumerate}
\item  OR $k+1$ lies  strictly to the left of $k$ and
   \begin{enumerate}
   \item  $k+2$ lies  weakly to the right of $k+1$ and strictly left of $k$, or
   \item  $k-1$ lies  strictly to the right of $k+1$ and weakly left of $k$.
   \end{enumerate}
\end{enumerate}

We proceed by induction on $|\alpha|$.
Clearly the proposition holds for $|\alpha|\leq 3$ since there is only one SCT of each composition shape in each of those cases.
Define $R(\alpha)$ be the set of all cells of the diagram which could be removed to obtain a covered shape in $\cL_C$.
We can also think of $R(\alpha)$ as the set of positions that a `1' entry could possibly appear in an SCT of shape $\alpha$.

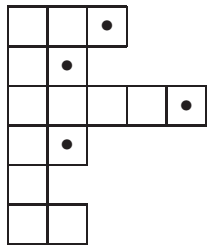
\begin{example}
\[
 \tableau{ \emt&\emt& \bullet \\ \emt & \bullet \\ \emt &\emt&\emt&\emt & \bullet \\ \emt & \bullet \\  \emt \\ \emt & \emt }
 \qquad\qquad R(3,2,5,2,1,2)= (3,2,5,2,1,2)\cskew(2,1,4,1,1,2)
\]
\end{example}

We can characterize $R(\alpha)$ as follows.
Each element of $R(\alpha)$ is a cell at the end of some row.
A cell in the first column of the diagram is in $R(\alpha)$ if and only if it is in position $(1,1)$ and row 1 has length 1.
A cell in column $j>1$ is in $R(\alpha)$ if and only if it is at the end of its row and it has no `hole' above it in the same column, that is, no position where a cell could be added to obtain a cover of $\alpha$.
Given $s\in R(\alpha)$ we write $\alpha\setnot s$ to denote the diagram (shape) obtained by removing $s$.
Similarly, we write $\alpha\cup s$ to denote the diagram obtained by adding $s$ to $\alpha$.

As earlier, given  $T\in SCT(\alpha)$, viewed as $T:\alpha\to[n]$,
and an integer $k$, we write $T+k$ to denote the $k$-shift of $T$ given by 
adding $k$ to each entry of the tableau.
Now given $s\in R(\alpha)$, consider the set of composition tableaux
\[ \cA(\alpha,s) := \{ T\in SCT(\alpha) \suchthat T(s)=1 \}.  \]
We note that $\cA(\alpha,s)$ can be constructed from the tableaux with one fewer cell, namely
\[ \cA(\alpha,s) = \{ (T+1)\cup(s\mapsto 1) \suchthat  T\in SCT(\alpha\setnot s) \}.  \]
By our induction hypothesis $G_{\alpha\setnot s}$ is connected.
Furthermore,  the dual Knuth move that exchanges $k$ and $k+1$ in $T$ is equivalent to the move that exchanges $k+1$ and $k+2$ in $T+1$.
It follows that the subgraph of $G_\alpha$ induced by $\cA(\alpha,s)$ (which, abusively, we also refer to as $G_{\alpha\setnot s}$) is connected.
If $|R(\alpha)|=1$ then we are done.

Otherwise, it remains to show that these $G_{\alpha\setnot s}$ subgraphs of $G_\alpha$ are connected to each other.
Consider a derived graph $\widehat{G}$ with vertex set $R(\alpha)$
where $s_1$ and $s_2$ are connected by an edge in $\widehat{G}$ if and only if  $G_{\alpha\setnot s_1}$ and $G_{\alpha\setnot s_2}$ are directly connected in $G_\alpha$.
We claim that $\widehat{G}$ is connected.
Let us partition $R(\alpha)$ into two sets $R(\alpha)=X\cup Y$ where $X$ contains every cell in $R(\alpha)$ that is the highest element of $R(\alpha)$ in its respective column, and $Y$ contains the remaining cells.
We will show that the subgraph of $\widehat{G}$ induced by $X$ is complete (a clique), and that every vertex in $Y$ is connected to a vertex in $X$, thus showing that $\widehat{G}$ is connected.
To show that $G_{\alpha\setnot s_1}$ and $G_{\alpha\setnot s_2}$ are directly connected in $G_\alpha$, that is, that  $s_1,s_2$ are connected in $\widehat{G}$,
it suffices to show that there exist $T_1\in\cA(\alpha,s_1)$ and $T_2\in\cA(\alpha,s_2)$ such that $T_1$ and $T_2$ differ by an elementary dual Knuth move.

\emph{Case: $s_1, s_2\in X$  (different columns).}
Without loss of generality we assume that $s_1$ is to the left of $s_2$ in the diagram.
Let $s_3$ be the cell in the same row as $s_2$ and to its immediate left.
Let $T\in SCT(\alpha\setnot\{s_1,s_2,s_3\})$ and set $T_1=(T+3)\cup(s_i\mapsto i)$ and let $T_2$ be the result of the dual Knuth move that exchanges  entries 1 and 2.
The column words of the tableaux are of the form $\cdots1\cdots 3 \cdots 2 \cdots$ and  $\cdots2\cdots 3 \cdots 1 \cdots$ respectively.
Examples:

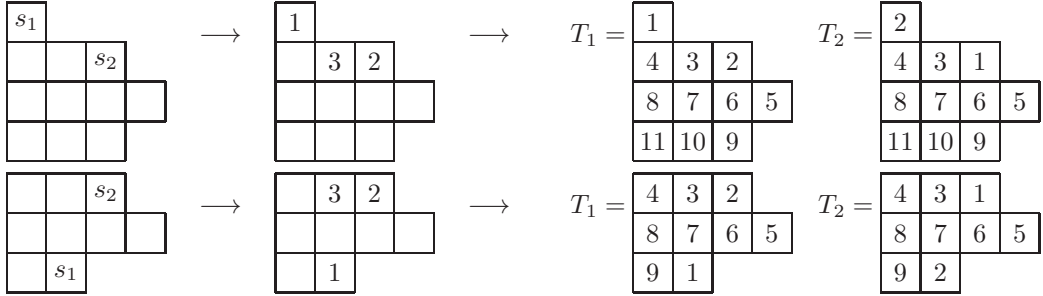
\[ \tableau{ s_1 \\ \emt & \emt & s_2 \\ \emt &\emt &\emt &\emt \\ \emt &\emt &\emt }
\quad \longrightarrow \quad
\tableau{ 1 \\ \emt & 3 & 2 \\ \emt &\emt &\emt &\emt \\ \emt &\emt &\emt }
\quad \longrightarrow \qquad
T_1 = \tableau{ 1 \\ 4 & 3 & 2 \\ 8 &7 &6 &5 \\ 11 &10 &9 }
\quad T_2 = \tableau{ 2 \\ 4 & 3 & 1 \\ 8 &7 &6 &5 \\ 11 &10 &9 }
\]

\[ \tableau{  \emt & \emt & s_2 \\ \emt &\emt &\emt &\emt \\ \emt & s_1  }
\quad \longrightarrow \quad
\tableau{  \emt & 3 & 2 \\ \emt &\emt &\emt &\emt \\ \emt &1 }
\quad \longrightarrow \qquad
T_1 = \tableau{  4 & 3 & 2 \\ 8 &7 &6 &5 \\ 9 &1 }
\quad T_2 = \tableau{  4 & 3 & 1 \\ 8 &7 &6 &5 \\ 9 &2 }
\]
\medskip
This case demonstrates that the subgraph of $\widehat{G}$ induced by $X$ is complete.

\medskip
\emph{Case: $s_1\in X, s_2\in Y$, both in column $j$.}
Clearly $j>1$.
Let $s_3$ be the cell in the same row as $s_1$ and to its immediate left.
Let $T\in SCT(\alpha\setnot\{s_1,s_2,s_3\})$
and set \[ T_1=(T+3)\cup(s_1 \mapsto 1, s_2 \mapsto 3, s_3 \mapsto 2)\]  and let $T_2$ be the result of the dual Knuth move that exchanges entries 3 and 2, which in fact rotates the three entries 3, 2, and 1 in the SCT.
The column words  of the tableaux are of the form $\cdots2\cdots 13 \cdots$ and  $\cdots3\cdots 12 \cdots$ respectively.
Example:

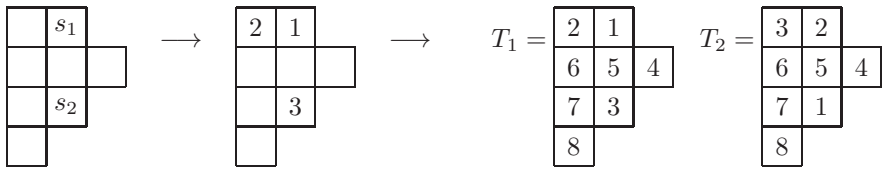
\[ \tableau{  \emt & s_1  \\ \emt &\emt &\emt    \\ \emt & s_2 \\ \emt  }
\quad \longrightarrow \quad
\tableau{ 2 & 1  \\ \emt &\emt &\emt    \\ \emt & 3 \\ \emt  }
\quad \longrightarrow \qquad
T_1 = \tableau{  2 & 1  \\ 6 &5 &4    \\ 7 & 3 \\ 8  }
\quad T_2 = \tableau{ 3 & 2  \\ 6 &5 &4   \\ 7 & 1 \\ 8  }
\]

\medskip
This case demonstrates that every vertex $s\in Y$ is connected to a vertex in $X$. 
As above, these two cases establish connectivity of $\widehat{G}$, and hence of $G_\alpha$.
\end{proof} 

%
%
\section{Applications of skew quasisymmetric Schur functions}\label{sec:apps}

\subsection{Symmetric skew quasisymmetric Schur functions} \label{subsec:sym-SQS}

Some skew quasisymmetric Schur functions are symmetric.
For example, if you take a skew SSRT and extend the base shape by adding an extra column of cells on the left, one for every row, the resulting skew SSRT also meets the definition of an SSCT.
Conversely, an SSCT of shape $\gamma\cskew\mu$ where $\mu$ is a partition and $\ell(\mu)=\ell(\gamma)$ will have strictly decreasing column entries, and thus be an SSRT as well, and the tableau obtained by removing the first column of cells from its base shape will still be an SSRT.
The combinatorial formulas  \eqref{eqn:skewSchurComb}  and \eqref{eqn:qschur-stdtabs-skew} then imply that every skew Schur function is equal to a skew quasisymmetric Schur function.
\begin{example}
\[ \begin{array}{c}
\tableau{\gst&\gst & 3 \\ \gst &4 & 1 \\ 2} \qquad\longleftrightarrow\qquad
  \tableau{ \gst & \gst&\gst & 3 \\   \gst &\gst &4 & 1 \\\gst & 2  }
\vspace{3pt} \\  s_{(3,3,1)/(2,1)} \;=\; \qs_{(4,4,2)\cskew(3,2,1)}
\end{array} \]
\end{example}

More generally, say that a skew composition shape $\gamma\cskew\beta$ is \emph{uniform} if all of the rows of the skew shape that have a cell in the first column are of the same length, that is, $\gamma_i=\gamma_j$ for all $1\leq i<j\leq \ell(\gamma)-\ell(\beta)$.
For example, the shape $(3,3,3,6,2,3)\cskew(2,1,1)$ is uniform.
A consequence of the proof of Proposition \ref{prop:C-complete} is that, since SCT of uniform shape cannot have any rigid row pair, all dual Knuth moves applied to an SCT of uniform shape result in another SCT of the same shape.
It follows that the SCT of that shape, or rather their SRT images under the bijection $\rho$, can be partitioned into complete dual equivalence classes, in Haiman's sense \cite{haiman-dualeq}.
Hence, we have the following.
\begin{corollary} \label{prop:sym-SQS}
Let $\gamma\cskew\beta$ be a uniform skew composition shape.
Then $\qs_{\gamma\cskew\beta}$ is symmetric, and expands as a nonnegative integer linear combination of Schur functions.
\end{corollary}

\subsection{The algebra of Poirier-Reutenauer and free Schur functions}\label{subsec:PR}

Poirier and Reute\-nauer \cite{poirier-reutenauer} introduced a dual pair of noncommutative Hopf algebras whose bases are 
parameterized by straight SYT.
Of these, the one we consider here, which we designate $PR$,  has been shown to be isomorphic to $FSym$, the algebra of \emph{free Schur functions} defined by Duchamp, Hivert, and Thibon \cite{ncsf6}.
\hilight{In \cite[Theorem 4.3]{poirier-reutenauer},} 
 it is shown that $\Sym$ is a quotient of $PR$, the linear map being determined by  
\[ T\mapsto s_\lambda, \qquad \text{ where } sh(T)=\lambda. \]
We show that this morphism of algebras factors through $\Nsym$.
\hilight{We emphasize that while the map $PR\to\Sym$ of Poirier and Reutenauer and the forgetful map $\chi:\Nsym\to\Sym$ are both Hopf algebra morphisms, the map we present below is only an anti-morphism of algebras since it does not respect the coproduct.}

Here we use straight SRT for basis 
elements.
The product of basis elements $T_1$ and $T_2$ in $PR$ is defined by the shifted shuffle of the Knuth equivalence classes of permutations that they index.
The effect is that if $T_2$ has $n$ cells and $sh(T_1)=\mu$, their product is
\[ T_1* T_2 \;=\; \sum_T T, \]
where the sum runs over all SRT $T$ such that $T|_\mu = T_1+n$ and $rect(T|_{\nu/\mu})=T_2$, where $sh(T)=\nu$.
That is, $T$ restricted to the base shape $\mu$ is $T_1+n$ and the rectification of the remaining skew tableau, as an SRT, is $T_2$.
\begin{example}
\[ \tableau{3&2\\1} * \tableau{3&2&1} \;=\;
\tableau{6&5 &3&2&1\\ 4} \,+\, \tableau{6&5 &2&1\\ 4&3} \,+\, \tableau{6&5 &2&1\\ 4\\3} \,+\, \tableau{6&5 &1\\ 4&2\\3} \]
\end{example}

\begin{theorem} \label{prop:PR}
The linear map  
 $\varphi:PR\to \Nsym$ given by \[ \varphi(T) \;=\; \nqs_\alpha, \qquad \text{ where } \alpha =\cshape(T), \]
is a surjective anti-morphism of algebras, i.e., \begin{equation} \label{eqn:PR}
 \varphi(T_1*T_2) \;=\; \varphi(T_2)\cdot \varphi(T_1).
\end{equation}
\end{theorem}
\begin{proof}
Suppose $T_1,T_2\in SRT$. 
As described above, in $PR$, the product $T_1*T_2$ is a positive sum of SRT, hence $\varphi(T_1*T_2)$ is a positive sum of noncommutative Schur functions. 
By Theorem \ref{thm:qschur-skew-pos}, the right hand side $\varphi(T_2)\cdot \varphi(T_1)$ is also a positive sum of noncommutative Schur functions.
We show that (1) there is a $C$-shape preserving bijection between the terms of $T_1*T_2$ and the terms of $\varphi(T_2)\cdot \varphi(T_1)$, 
and (2) the map $\varphi$ respects products, i.e., $\varphi(T_1*T_2)$ depends only on $\cshape(T_1)$ and $\cshape(T_2)$.

Suppose that $\cshape(T_1)=\beta$,  $| \beta |=m$, $\cshape(T_2)=\alpha$, and $| \alpha |=n$.
Let $T$ be a term of $T_1*T_2$, say with $\cshape(T)=\gamma$.
From the construction in the proof of Proposition \ref{prop:rho-general} using the notation in the proof of Proposition~\ref{prop:qschur-skew}, it is straightforward to show that  $\rho^{-1}(T)|_\beta-n = \Omega_m(\rho^{-1}(T))= \rho^{-1}(T_1)$ and $rect(\rho^{-1}(T)|_{\gamma\cskew\beta})=rect(\mho_n(\rho^{-1}(T)))=\rho^{-1}(T_2)$.
So for each term $T$, we consider the $C$-equivalence class of SCT $[\hat{T}]_C$ where $\hat{T}=\mho_n(\rho^{-1}(T))$. 
Note that the rectified composition shape of $[\hat{T}]_C$ is $\alpha$  since $rect(\hat{T})=\rho^{-1}(T_2)$.
By  Proposition  \ref{prop:C-complete}, each $C$-equivalence class of SCT of rectified shape $\alpha$ contains exactly one member that rectifies to $\rho^{-1}(T_2)$, so these equivalence classes $[\hat{T}]_C$ are distinct across terms of $T_1*T_2$.

Also by Proposition \ref{prop:C-complete}, each of these equivalence classes $[\hat{T}]_C$ contains exactly one member that rectifies to $U_\alpha$.
Conversely, suppose that $S$ is an SCT of shape $\eta\cskew\beta$ for some composition $\eta$ such that $rect(S)=U_\alpha$.
Again by Proposition \ref{prop:C-complete}, $[S]_C$ contains a unique member $\check{S}$ that rectifies to $\rho^{-1}(T_2)$, hence $\rho(\rho^{-1}(T_1+n)\cup\check{S})$ is one of the terms of $T_1*T_2$.
Thus the terms of $T_1*T_2$ are in bijection with the set 
\[ \{ S\in SCT : sh(S)=\gamma\cskew\beta \text{ for some } \gamma, \text{ and } rect(S)=U_\alpha \} \]
and the bijection preserves the overall composition shape $\gamma$ of the terms.
By Theorem \ref{thm:qschur-skew-pos}, the terms of $\nqs_\alpha\cdot \nqs_\beta$, that is, the expansion of the right hand side of \eqref{eqn:PR}, when considered as a sum of terms each with coefficient 1,  are also in bijection with this set, the bijection preserving the overall shape $\gamma$ of each term.
This establishes the desired bijection.

Consideration of the above bijection shows that it only depends on $\cshape(T_1)$ and $\cshape(T_2)$, and not on the specific tableaux (fillings) of those shapes.
Thus the map $\varphi$ respects products as desired.
\end{proof}
\hilight{}

\subsection{$NCSym$, $NCQSym$ and their noncommutative Schur functions}
\label{subsec:NCSym}

In this subsection we extend the definition of Schur functions in noncommuting variables studied by 
Rosas and Sagan in \cite{rosas-sagan} to quasisymmetric Schur functions in
noncommuting variables, and then prove that they project naturally onto quasisymmetric Schur functions. For this we need to consider two subalgebras of $\mathbb{Q}\ll x_1, x_2, \ldots \gg$, the Hopf algebra of formal power series in noncommuting variables.

For the first subalgebra, let $[n]=\{1, 2, \ldots ,n\}$. Then a \emph{set partition} of $[n]$ is a family $\pi = \{ A_1, A_2, \ldots , A_\ell\}$ of pairwise disjoint nonempty sets such that
$\cup _{i=1} ^{\ell} A _i=[n]$, and is denoted by $\pi \vdash [n]$. If $|A_i| = \alpha _i$ then let $\lambda (\pi)$ denote the partition of $n$ determined by $\alpha _1, \ldots , \alpha _\ell$. Given a set partition, $\pi = \{ A_1, A_2, \ldots , A_\ell\} \vdash [n]$, define the \emph{monomial symmetric function in noncommuting variables} $\ncm _\pi$ to be

$$\ncm _\pi = \sum _{(i_1, \ldots , i_n)} x_{i_1}\cdots x_{i_n},$$
where $i_j = i_k$ if and only if $j, k \in A_m$ for some $1\leq m \leq \ell$.

\begin{example} If $n=3$ and $\pi = \{13, 2\}$, then $\lambda(\pi)=(2,1)$ and
$$\ncm _\pi = x_1x_2x_1 + x_2x_1x_2 + x_1x_3x_1 + x_3x_1x_3 +  \cdots .$$
\end{example}
The Hopf algebra $NCSym$ is then defined as
$$NCSym = \bigoplus_{n\geq 0} \spam \{ \ncm _\pi \ : \ \pi \vdash [n] \}$$
and its structure has been
 studied in \cite{bergeron-zabrocki-2, bergeron-zabrocki-1, rosas-sagan}.
Let a \emph{dotted reverse tableau} $\Tdot$ of shape $sh (\Tdot)$ be an SSRT $T$ of shape $sh(T)$ 
which has for each $k=1, \ldots , | sh (T) |$ 
exactly one entry with $k$ dots placed above it.  
Then \cite{rosas-sagan} defined the \emph{Schur function in noncommuting variables} $\SRS _\lambda $ to be
$$\SRS_\lambda = \sum _{sh(\Tdot) = \lambda} x^{\Tdot},$$
where the sum is over all dotted reverse tableaux $\Tdot$ of shape $\lambda$, and $x^{\Tdot}$ is the monomial with $x_i$ in position $j$ if and only if $\Tdot$ has a cell containing $i$ with $j$ dots above it.

\begin{example} Restricting ourselves to 2 variables,
$$\SRS _{21} = 2 x_2x_1x_1  + 2 x_1x_2x_1+2x_1x_1x_2 +  2 x_1x_2x_2 + 2 x_2x_1x_2 + 2x_2x_2x_1$$from the dotted reverse tableaux
$$\tableau{\dot{2}&\ddot{1}\\\dddot{1}}\ \tableau{\dot{2}&\dddot{1}\\\ddot{1}}\  \tableau{\ddot{2}&\dot{1}\\\dddot{1}}\  \tableau{\ddot{2}&\dddot{1}\\\dot{1}}\  \tableau{\dddot{2}&\dot{1}\\\ddot{1}}\  \tableau{\dddot{2}&\ddot{1}\\\dot{1}}\
\tableau{\ddot{2}&\dddot{2}\\\dot{1}} \    
\tableau{\dddot{2}&\ddot{2}\\\dot{1}} \   
\tableau{\dot{2}&\dddot{2}\\\ddot{1}} \   
\tableau{\dddot{2}&\dot{2}\\\ddot{1}} \   
\tableau{\dot{2}&\ddot{2}\\\dddot{1}} \   
\tableau{\ddot{2}&\dot{2}\\\dddot{1}}.$$
\end{example}
Considering the {forgetful map}
$\chi : \mathbb{Q} \ll x_1, x_2, \ldots \gg\ \rightarrow   \mathbb{Q} [[ x_1, x_2, \ldots ]]$, which lets the variables commute, 
we have 
$\chi (\SRS _\lambda )= n! s_\lambda$ where $|\lambda | = n$ and furthermore $\{ \SRS _\lambda \} _{\lambda \vdash n\geq 0} \subseteq NCSym$ by expressing $\SRS _\lambda$ in terms of the $\ncm _\pi$ combinatorially, 
see \cite{rosas-sagan} for details.

For the second subalgebra, a \emph{set composition} of $[n]$ is an \emph{ordered} family $\Pi = ( A_1, A_2, \ldots , A_\ell\ )$ of pairwise disjoint nonempty sets such that
$\cup _{i=1} ^{\ell} A _i=[n]$,  and is denoted by $\Pi \vDash [n]$. If $|A_i| = \alpha _i$ then let $\alpha (\Pi)$ denote the composition  $(\alpha _1,\cdots ,\alpha _\ell)\vDash n$.  Given a set composition, $\Pi = (A_1, A_2, \ldots , A_\ell) \vDash [n]$, define the \emph{monomial quasisymmetric function in noncommuting variables} $\ncM _\Pi$ to be

$$\ncM _\Pi = \sum _{(i_1, \ldots , i_n)} x_{i_1}\cdots x_{i_n},$$
where \begin{itemize}\item $i_j = i_k$ if and only if $j, k \in A_m$ for some $1\leq m \leq \ell$, and \item $i_j < i_k$ if and only if $j\in A_{m_1}, k\in A_{m_2}$ and $m_1 < m_2$.\end{itemize}

\begin{example} If $n=3$ and $\Pi = (2, 13)$, then $\alpha(\Pi)=(1,2)$ and
$$\ncM _\Pi =  x_2x_1x_2 +  x_3x_1x_3 +  \cdots .$$
\end{example}

The Hopf algebra $NCQSym$, introduced by 
\hilight{Aguiar and Mahajan \cite[Section 6.2.5]{aguiar-mahajan-1} and studied further by}
Bergeron and Zabrocki in \cite{bergeron-zabrocki-1},
is then defined as
$$NCQSym = \bigoplus _{n\geq 0} \spam \{ \ncM _\Pi \ : \ \Pi \vDash [n] \}\:.$$
In analogy to Schur functions in noncommuting variables, we can also define quasisymmetric Schur functions in noncommuting variables.

\begin{definition}\label{def:schurinnc} Let $\alpha$ be a composition, and let a \emph{dotted composition tableau} $\Tdot$ of shape $sh (\Tdot)$ be an SSCT $T$ of shape $sh(T)$ with exactly one entry with $k$ dots placed above it for $k=1, \ldots , | sh(T) |$. Then  the \emph{quasisymmetric Schur function in noncommuting variables} $\QSRS _\alpha $ 
is 
defined to be
$$\QSRS_\alpha = \sum _{sh(\Tdot) = \alpha} x^{\Tdot},$$ where the sum is over  all dotted composition tableaux $\Tdot$ of shape $\alpha$, and $x^{\Tdot}$ is the monomial with $x_i$ in position $j$ if and only if $\Tdot$ has a cell containing $i$ with $j$ dots above it.
\end{definition}

\begin{example} Restricting ourselves to 2 variables,
$$\QSRS _{12} = 2 x_1x_2x_2 + 2 x_2x_1x_2 + 2x_2x_2x_1$$from the dotted composition tableaux
$$
\tableau{\dot{1}\\\ddot{2}&\dddot{2}}\
\tableau{\dot{1}\\\dddot{2}&\ddot{2}}\
\tableau{\ddot{1}\\\dot{2}&\dddot{2}}\
\tableau{\ddot{1}\\\dddot{2}&\dot{2}}\
\tableau{\dddot{1}\\\dot{2}&\ddot{2}}\
\tableau{\dddot{1}\\\ddot{2}&\dot{2}}\ .$$
\end{example}
We observe that for any given composition tableau of weight $n$, there are clearly 
$n!$ dotted composition tableaux of this type.

Many results proved in \cite{rosas-sagan} for Schur functions in noncommuting variables can be extended to quasisymmetric Schur functions in noncommuting variables. For a composition $\alpha = (\alpha _1, \alpha _2, \ldots )$ let $\alpha ! = \alpha _1 ! \alpha _2 ! \cdots$.

\begin{theorem}\label{the:qsncom} For a partition $\lambda$ and compositions $\alpha, \beta$ where $n=|\lambda | = |\alpha | = |\beta|$,
\begin{enumerate}
\item $\SRS _\lambda = \sum _{\tilde{\alpha} = \lambda} \QSRS _\alpha$.
\item $\QSRS _\alpha = \sum _\beta \beta ! K_{\alpha\beta} \sum _{\alpha(\Pi) = \beta} \ncM _\Pi$, where  $K_{\alpha\beta}=$ the number of SSCT $T$ such that $sh(T)=\alpha$ and $cont(T)=\beta$.
\item The $\QSRS _\alpha$ are linearly independent.
\item $\chi(\QSRS _\alpha) = n! \qs _\alpha$ and $\tilde{\chi}(n ! \qs _\alpha) = \QSRS _\alpha$, where $\chi (\ncM _\Pi) = M _{\alpha (\Pi)}$ with 
\hilight{right} 
 inverse   $\tilde{\chi} (M _\alpha) = \sum _{\alpha(\Pi) = \alpha} \frac{\alpha !}{n!} \ncM _\Pi$.
\end{enumerate}
\end{theorem}

\begin{proof}
\begin{enumerate}
\item $$\SRS _\lambda = \sum _{sh(\rho(\Tdot))=\lambda} x^{\Tdot} =  \sum _{\tilde{\alpha} = \lambda} \sum _{sh(\Tdot)=\alpha} x^{\Tdot} =  \sum _{\tilde{\alpha} = \lambda} \QSRS _\alpha,$$
where the first sum is over all dotted composition tableaux $\Tdot$ that map naturally under $\rho$ to a dotted reverse tableau of shape $\lambda$.
\item Consider a monomial $x^{\Tdot}$ where $\Tdot$ is a dotted composition tableau with $sh(\Tdot) = \alpha$ appearing in $\ncM _\Pi$ where $\alpha (\Pi) = \beta$. The number of composition tableaux $T$ with $sh(T)=\alpha$ and $cont(T)=\beta$ is $K_{\alpha \beta}$. Since the number of ways to distribute dots to yield a dotted composition tableau with associated monomial $x^{\Tdot}$ is $\beta !$ the result follows.
\item This follows immediately from the previous part and \cite[Proposition 6.7]{HLMvW-1}.
\item It is straightforward to check, using a proof analogous to \cite[Proposition 4.3]{bergeron-zabrocki-2}, that the forgetful map $\chi$ is a surjective Hopf morphism 
    $\chi: NCQSym \rightarrow QSym$ that satisfies
$\chi (\ncM _\Pi) = M _{\alpha (\Pi)}$. 
The first equation follows immediately from the observation made before the 
theorem. 
It is also straightforward to check using a proof analogous to \cite[Proposition 4.1]{rosas-sagan} that 
$\tilde{\chi}: QSym \rightarrow NCQSym$ is an injective inclusion, and is a \hilight{right} inverse for $\chi$. Now from the second part, all $\ncM _\Pi$ with $\alpha(\Pi)=\beta$ have the same coefficient in $\QSRS _\alpha$ and so $\QSRS _\alpha$ is in the image of $\tilde{\chi}$. The result now follows from $\chi(\QSRS _\alpha) = n! \qs _\alpha$ 
as  $\tilde{\chi}$ is a \hilight{right} inverse for $\chi$.
\end{enumerate}
\end{proof}
\begin{remark} Theorem~\ref{the:qsncom} and \cite[Theorem 6.2]{rosas-sagan} yield the following commutative diagram, which relates the Schur functions of $Sym$, $QSym$, $NCSym$, and $NCQSym$.
$$\xymatrix{
\SRS _\lambda \ar@2{-}[r]  \ar@{->}[d]  _\chi& \sum _{\tilde{\alpha} = \lambda} \QSRS _\alpha \ar@{->}[d] ^\chi \\
n!\ s_\lambda \ar@2{-}[r] & n!\ \sum _{\tilde{\alpha} = \lambda} \qs _\alpha}$$
\end{remark}

%
\subsection{Pieri operators \hilight{and skew quasisymmetric Schur functions}}\label{subsec:pieri}

In \cite{bergeron-vW} the notion of a Pieri operator was introduced. More precisely, given a graded poset $P$ with rank function $rk: P \rightarrow \bP$ and $k\in \bP$ a \emph{Pieri operator} is a linear map
 $\po _k : \mathbb{Z}P \rightarrow \mathbb{Z}P$ such that for all $x\in P$ the support of $x.\po _k \in \mathbb{Z}P$ consists only of elements $y\in P$ such that $x<y$ and $rk(y) - rk(x) = k$. Furthermore they identified $\po _k$ with $\bh _k \in \Nsym$ and by duality established  a collection of homogeneous quasisymmetric functions associated to every interval $[ x, y ]$ of $P$
$$K_{[x,y]} = \sum _\alpha \langle x. \overline{a} _\alpha , y \rangle b _\alpha,$$
where $\langle \cdot \ ,\ \cdot \rangle$ is the bilinear form induced by the Kronecker delta function, $\{ a _\alpha \} _{\alpha \vDash n\geq 0}$ is a graded basis of $\Nsym$ and $\{ b _\alpha \}  _{\alpha \vDash n\geq 0}$ is the corresponding dual basis of $\Qsym$. Intuitively, the coefficient of $b _\alpha$ in $K _{[ x, y]}$ can be thought of as the number of saturated chains from $x$ to $y$ in $P$ satisfying conditions imposed by $\overline{a} _\alpha$. 
Depending on the choice of $P$ and $\po _k$ examples of $K_{[x,y]}$ include skew Schur functions, Stanley symmetric functions \cite{stanley-cox}, skew Schubert functions \cite{bergeron-sottile-2}, and the noncommutative Schur functions of Fomin and Greene \cite{fomin-greene-1}. 
To identify a further example, we need the descent Pieri operator arising in the following theorem.

\begin{theorem}\label{the:descentpieri}
\cite[Equation 4]{bergeron-sottile-1} Let $P$ be a graded edge labeled poset whose covers are labeled by elements of a totally ordered set $(B, <)$. Consider for $x\in P$ the \emph{descent Pieri operator}
$$x. \po _k  = \sum _{\omega} end(\omega),$$
where the sum is over all saturated chains $\omega$ of length $k$
$$\omega: x \overset{b_1}{\rightarrow} x _1\overset{b_2}{\rightarrow} \cdots  \overset{b_k}{\rightarrow} x_k=end(\omega)$$for $ b_1 \leq b_2 \leq \cdots \leq b_k \in B$.

Furthermore, given a saturated chain $\omega$ of length $n$ in $P$ with labels $b_1, b_2, \ldots , b_n \in B$ let its \emph{descent set} be $descents(\omega) = \{ i \ | \ b_i > b_{i+1}\}$, and let its corresponding \emph{descent composition} $Des(\omega)$  be the composition of $n$ defined by $set(Des(\omega)) = descents(\omega)$. Then
$$K_{[x,y]} = \sum _{\omega \in ch[x,y]} L_{Des(\omega)},$$ where $ch[x,y]$ is the set of all saturated chains from $x$ to $y$.
\end{theorem}

We can now identify a new example of $K_{[x,y]}$.

\begin{theorem}\label{the:pieriskewqs} 
Let $\cL _C '$ be the dual poset of $\cL_C $, i.e., $\beta$ covers $\gamma$ in $\cL _C '$ if and only if $\gamma$ covers $\beta $ in $\cL _C$.
Let $P$ be the poset $\mathcal{L} _C '$ with edges labeled
$$x \overset{(-i, -j)}{\longrightarrow} \tilde{x},$$
where $i$ (respectively, $j$) is the column (respectively, row) index of the cell, using cartesian coordinates, appearing in $x$ but not $\tilde{x}$.
Let these labels be totally ordered: $(i,j)<(k,\ell)$ if and only if $i<k$ or ($i=k=-1$ and $j>\ell$) or ($i=k<-1$ and $j<\ell$).

Then considering the descent Pieri operator on $P$ we have
$$K_{[\gamma, \beta]} = \qs _{\gamma \cskew\beta}.$$
\end{theorem}

\begin{proof} By Proposition \ref{prop:qschur-skew} we know
$$\qs _{\gamma \cskew \beta} = \sum _{T\in SCT (\gamma\cskew\beta)}L _{Des(T)}.$$
Therefore, by Theorem~\ref{the:descentpieri} it suffices to show
$$\sum _{\omega\in ch[\gamma ,\beta]}L _{Des(\omega)}=\sum _{T\in SCT (\gamma \cskew\beta)}L _{Des(T)},$$
which we do via a bijection between the  chains in $ch[\gamma , \beta]$ and  $SCT(\gamma \cskew\beta)$ that preserves descent sets of chains and tableaux.
By Proposition \ref{prop:sct-chains}, given  $\omega \in ch [\gamma, \beta]$ $$\omega: \gamma \overset{(-i _1, -j_1)}{\longrightarrow}  \cdots  \overset{(-i _\ell, -j_\ell)}{\longrightarrow} \beta,$$
there exists a corresponding $T^\omega \in SCT (\gamma \cskew\beta)$ with $k$ in cell $(j_k,i_k)$ for all $1\leq k \leq \ell$.

Finally, we need to check that $(-i _k, -j_k) \in descents(\omega) \Leftrightarrow k\in descents(T^\omega)$. Note that $(-i _k, -j_k) \in descents(\omega)$ if and only if
$k+1$ is in a column weakly to the right of $k$ in $T^\omega$ and hence by definition $k\in descents(T^\omega)$.
\end{proof}

\hilight{\begin{remark}
Theorem \ref{the:pieriskewqs}  could also be established using the universal property of $\Qsym$ described in \cite{aguiar-bergeron-sottile}.
\end{remark}}

\begin{remark} The aforementioned noncommutative Schur functions of Fomin and Greene \cite{fomin-greene-1} give rise to symmetric functions $F _{y/x}$ and these in turn are another example of $K_{[x,y]}$ arising from descent Pieri operators, this time with underlying labeled multigraph $x \overset{-i}{\rightarrow} {x.u_i}$. For further details see \cite[Example 6.4]{bergeron-vW}.
\end{remark}

\section{Further avenues} \label{sec:conc}

The noncommutative Littlewood-Richardson rule in Theorem  \ref{thm:qschur-skew-pos}, in addition to the quasisymmetric Littlewood-Richardson rule presented in \cite{HLMvW-2} and the quasisymmetric Kostka numbers identified in \cite{HLMvW-1}, raises the question of what other classical Schur function properties lift to quasisymmetric or noncommutative Schur functions. For example, can the Jacobi-Trudi determinant formula for computing skew Schur functions be generalized to quasisymmetric skew Schur functions, or can a determinantal formula be found for noncommutative Schur functions using quasideterminants that arise in the study of $NSym$? 
Another example of a question to pursue is, since Schur functions arise naturally as irreducible characters in the representation theory of the symmetric group, 
whether representation 
theoretic interpretations exist for either quasisymmetric or nonsymmetric Schur functions.
\hilight{Certainly,  a representation theoretic interpretation of $\Qsym$ exists which involves the $0$-Hecke algebra, via fundamental quasisymmetric functions \cite{ncsf4}.  
Since quasisymmetric Schur functions are nonnegative linear combinations of fundamental quasisymmetric functions \cite[Proposition 5.2]{HLMvW-1}, quasisymmetric Schur functions would correspond to certain representations of the $0$-Hecke algebra, and it would be interesting to know precisely which ones.}

Closely related to quasisymmetric Schur functions are  refinements of them known as Demazure atoms, and also related are Demazure characters that consist of linear combinations of Demazure atoms and arise in the study of Schubert calculus and other areas. In \cite{HLMvW-2} it was shown that a Schur function multiplied by a quasisymmetric Schur function, Demazure atom, or Demazure character, and expanded in the same basis exhibited a refined Littlewood-Richardson rule. Therefore, due to the similarities between quasisymmetric Schur functions, Demazure atoms, and Demazure characters, another avenue to pursue is  properties of skew Demazure atoms or characters, and then skew Macdonald polynomials. The latter polynomials would arise through the symmetrization of Demazure atoms and the introduction of additional parameters $q,t$.

Considering symmetrization, we can also pursue the classification of when a skew quasisymmetric Schur function $\qs _{\gamma\cskew\beta}$ is symmetric.
In this regard, we  conjecture the converse of Corollary \ref{prop:sym-SQS}.

\begin{conjecture}
Suppose  $\qs_{\gamma\cskew\beta}$ is symmetric.
Then $\gamma\cskew\beta$ is a uniform skew composition shape.
\end{conjecture}

Enlarging the scope of our questions we can ask what properties are possessed  by $\mathcal{L} _C$. Although Remark  \ref{rem:LC} observes $\cL_C$ is not a lattice, it may possess other interesting properties.
Meanwhile, turning our attention to $NCSym$ and $NCQSym$ we can investigate whether there exist refinements of (quasisymmetric) Schur functions in noncommuting variables that form a basis for $NCSym$ and $NCQSym$. Lastly, returning to the diagram in Section~\ref{sec:intro}, we can explore maps to algebras related to these, and discover where these maps take quasisymmetric and nonsymmetric Schur functions.


\def\cprime{$'$}

\end{document}